%%%%%%%%%%%%%%%%%%%%%%%%%%%%%%%%%%%%%%%%%%%%%%%%%%%%%%%%%%%%%%%%%%%%%%%%%%
%                                                                        %
%                                                                        % 
%                  by J. Bouttier and  E. Guitter                        %
%                TEX file, using lanlmac.tex macros                      %
%                                                                        %
%                                                                        %
%%%%%%%%%%%%%%%%%%%%%%%%%%%%%%%%%%%%%%%%%%%%%%%%%%%%%%%%%%%%%%%%%%%%%%%%%%
\input lanlmac
%\draftmode
\def\href#1#2{{#2}}

\input epsf.tex

\overfullrule=0mm

\newcount\figno
\figno=0
\def\fig#1#2#3{
\par\begingroup\parindent=0pt\leftskip=1cm\rightskip=1cm\parindent=0pt
\baselineskip=11pt
\global\advance\figno by 1
\midinsert
\epsfxsize=#3
\centerline{\epsfbox{#2}}
\vskip 12pt
{\bf Fig.\ \the\figno:} #1\par
\endinsert\endgroup\par
}
\def\figlabel#1{\xdef#1{\the\figno}}
\def\encadremath#1{\vbox{\hrule\hbox{\vrule\kern8pt\vbox{\kern8pt
\hbox{$\displaystyle #1$}\kern8pt}
\kern8pt\vrule}\hrule}}

%Macros 
%%%%%%%%%%%%%%%%%%%%%%%%%%%%%%%%%%%%%%%%%%%%%%%%%%%%%%%%%%%%%%%%%

\magnification=\magstep1
\baselineskip=12pt
\hsize=6.3truein
\vsize=8.7truein
 at 8truept
 at 8truept
 at 10truept
%\footline={\footsc the electronic journal of combinatorics
%   {\footbf 11} (2004), \#R00\hfil\footrm\folio}

%%%%%%%%%%%%%%%%%%%%%%%%%%%%%%%%%%%%%%%%%%%%%%%%%%%%%%%%%%%%%%%%%%%%%%%%
\font\bigrm=cmr12 at 14pt \centerline{\bigrm Planar maps and continued 
fractions}

\bigskip\bigskip

\centerline{J. Bouttier and E. Guitter}
  \smallskip
  \centerline{Institut de Physique Th\'eorique}
  \centerline{CEA, IPhT, F-91191 Gif-sur-Yvette, France}
  \centerline{CNRS, URA 2306}
\centerline{\tt jeremie.bouttier@cea.fr}
\centerline{\tt emmanuel.guitter@cea.fr}

  \bigskip

     \bigskip\bigskip

     \centerline{\bf Abstract}
     \smallskip
     {\narrower\noindent

We present an unexpected connection between two map enumeration
problems. The first one consists in counting planar maps with a
boundary of prescribed length. The second one consists in counting
planar maps with two points at a prescribed distance. We show that,
in the general class of maps with controlled face degrees, the
solution for both problems is actually encoded into the same quantity,
respectively via its power series expansion and its continued fraction
expansion. We then use known techniques for tackling the first problem
in order to solve the second. This novel viewpoint provides a
constructive approach for computing the so-called distance-dependent
two-point function of general planar maps. We prove and extend some
previously predicted exact formulas, which we identify in terms of
particular Schur functions. 
\par}

     \bigskip

%references
\nref\TutteCPM{W.T. Tutte, {\it A Census of Planar Maps}, Canad. J. of Math. 
{\bf 15} (1963) 249-271.}
\nref\TUTEQ{W.T. Tutte, {\it On the enumeration of planar maps}, Bull.\
Amer.\ Math.\ Soc.\ {\bf 74} (1968) 64-74.}
\nref\DGZ{See for instance: P. Di Francesco, P. Ginsparg and 
J. Zinn--Justin, {\it 2D Gravity and Random Matrices},
Physics Reports {\bf 254} (1995) 1-131, and references therein.}
\nref\GouldJack{See for instance Sect.~2.9 of: I.P. Goulden and D.M.
Jackson, {\it Combinatorial Enumeration}, John Wiley \& Sons, New
York, 1983, republished by Dover, New York, 2004, and references
therein.}
\nref\BIPZ{E. Br\'ezin, C. Itzykson, G. Parisi and J.-B. Zuber, {\it Planar
Diagrams}, Comm. Math. Phys. {\bf 59} (1978) 35-51.}
\nref\Schae{G. Schaeffer, {\it Conjugaison d'arbres
et cartes combinatoires al\'eatoires}, PhD Thesis, Universit\'e 
Bordeaux I (1998).}
\nref\CENSUS{J. Bouttier, P. Di Francesco and E. Guitter, {\it Census of planar
maps: from the one-matrix model solution to a combinatorial proof},
Nucl. Phys. {\bf B645}[PM] (2002) 477-499, arXiv:cond-mat/0207682.}
\nref\MOB{J. Bouttier, P. Di Francesco and E. Guitter. {\it 
Planar maps as labeled mobiles},
Elec. Jour. of Combinatorics {\bf 11} (2004) R69, arXiv:math.CO/0405099.}
\nref\Bouttier{J. Bouttier, {\it Physique statistique des surfaces al\'eatoires
et combinatoire bijective des cartes planaires}, PhD Thesis, Universit\'e
Paris 6 (2005).}
\nref\CS{P. Chassaing and G. Schaeffer, {\it Random Planar Lattices and 
Integrated SuperBrownian Excursion}, 
Probability Theory and Related Fields {\bf 128(2)} (2004) 161-212, 
arXiv:math.CO/0205226.}
\nref\GEOD{J. Bouttier, P. Di Francesco and E. Guitter, {\it Geodesic
distance in planar graphs}, Nucl. Phys. {\bf B663}[FS] (2003) 535-567, 
arXiv:cond-mat/0303272.}
\nref\MierRev{See for instance: G. Miermont, {\it Random maps and their 
scaling limits}, in C. Bandt, P. M\"orters, M. Z\"ahle (Eds.),
Proceedings of the conference Fractal Geometry and Stochastics IV, 
Greifswald (2008), 
Progress in Probability, Vol.\ {\bf 61}, 197-224, Bikha\"user (2009), and
references therein.}
\nref\AW{J. Ambj\o rn and Y. Watabiki, {\it Scaling in quantum gravity},
Nucl.Phys. {\bf B445} (1995) 129-144, arXiv:hep-th/9501049.}
\nref\FLAJFRAC{P. Flajolet, {\it Combinatorial aspects of continued
fractions}, Discrete Mathematics 32 (1980), 125--161.
Reprinted in the 35th Special Anniversary Issue of Discrete
Mathematics, Vol.\ {\bf 306}, Issue 10-11, Pages 992-1021 (2006).}
\nref\FLAJFRACBIS{P. Flajolet and R. Sedgewick, {\it Analytic combinatorics},
Cambridge University Press, 2009, section V.4.}
\nref\Wall{H.S. Wall, {\it Analytic Theory of Continued Fractions},
D. Van Nostrand, New York, 1948, reprinted by Chelsea, New York,
1973.}
\nref\FULHA{W. Fulton and J. Harris, {\it Representation Theory}, 
Springer-Verlag, New York, 1991.}
\nref\FULKRA{M. Fulmek and C. Krattenthaler, {\it Lattice path proofs
for determinant formulas for symplectic and orthogonal characters},
J. Combin.\ Theory Ser.\ {\bf A 77} (1997), 3-50.}
\nref\Viennot{X.G. Viennot, {\it Une th\'eorie combinatoire des polyn\^omes 
orthogonaux}, Lecture Notes UQAM, 217p., Publication du LACIM, Universit\'e 
du Qu\'ebec \`a Montr\'eal (1984), re\'ed. 1991.}
\nref\PDFRaman{P. Di Francesco, {\it Geodesic Distance in Planar Graphs: 
An Integrable Approach}, The Ramanujan Journal {\bf 10} (2005) 153-186,
arXiv:math/0506543 [math.CO].}
\nref\CMS{G. Chapuy, M. Marcus and G. Schaeffer, 
{\it A bijection for rooted maps on orientable surfaces}, 
SIAM J. Discrete Math. {\bf 23}(3) (2009) 1587-1611, 
arXiv:0712.3649 [math.CO].}
\nref\BECA{E.A. Bender and E.R. Canfield, {\it The number of
  degree-restricted rooted maps on the sphere}, SIAM J.\ Discrete
  Math.\ {\bf 7} (1994) 9-15.}
\nref\BMJ{M. Bousquet-M\'elou and A. Jehanne,
  {\it Polynomial equations with one catalytic variable, algebraic
    series and map enumeration}, J. Combin.\ Theory
  Ser.\ {\bf B 96} (2006) 623-672.}
\nref\Bro{W.G. Brown, {\it On the existence of square roots in certain
 rings of power series}, Math. Ann. {\bf 158} (1965) 82-89.}
\nref\DFG{P. Di Francesco and E. Guitter, {\it Integrability of graph 
combinatorics via random walks and heaps of dimers}, 
J. Stat. Mech. (2005) P09001, arXiv:math/0506542 [math.CO].}
\nref\LGV{See for instance: I.M. Gessel and X.G. Viennot, 
{\it Binomial determinants, paths and hook length formulae}, Adv. in Math. 
{\bf 58} (1985) 300-321; 
I.M. Gessel and X.G. Viennot, {\it Determinants, Paths, and
Plane Partitions}, preprint (1989), available for instance
at http://people.brandeis.edu/$\sim$gessel/.}
\nref\PSEUDOQUAD{J. Bouttier and E. Guitter, {\it Distance statistics
in quadrangulations with a boundary, or with a self-avoiding loop},
J. Phys.\ A: Math.\ Theor.\ {\bf 42} (2009) 465208, arXiv:0906.4892
[math-ph].}

%text
\newsec{Introduction}

\subsec{General introduction}

Maps are fundamental objects of combinatorial theory, introduced by
Tutte in the 60's [\xref\TutteCPM,\xref\TUTEQ], which also appeared later 
as natural models for random surfaces in physics \DGZ. 
Many questions about maps boil down to
enumeration problems, which in turn received a lot of attention from
different communities. Different techniques of enumeration have been
developed, for instance combinatorialists' quadratic method \GouldJack,
physicists'
matrix integrals \BIPZ\ and more recently bijective coding by trees 
[\xref\Schae-\xref\Bouttier]. 
Most results deal with global enumeration problems, say counting families
of maps with a control on their size and topology. The coding by trees
allowed to address more refined questions about the distance in maps,
say counting maps with a prescribed radius \CS\ or with marked pointed at
prescribed distances \GEOD. This approach led to important applications in
probability theory for the rigorous construction of scaling limits
of large random maps \MierRev.

One of the prominent global quantities, considered both in the recursive 
decomposition and the matrix integral approach, is the number of maps
with one boundary, with a control on both the total map size and the
boundary length. More precisely, this information is captured by the 
{\it moments}, which are generating functions of maps with a fixed boundary
length. These moments are determined by a closed system of equations, 
called either Tutte's or loop equations respectively in the combinatorial 
or physical communities.

On the other hand, one of the simplest observables related to distance
is the so-called {\it two-point function}, which is the 
number of maps with given size and having two marked vertices at fixed
distance. In particular the dependence of this two-point function over
the distance gives directly the average profile of random maps. The
distance-dependent two-point function was first considered in \AW\ in
the case of planar triangulations, where its scaling form was
predicted via a ``transfer-matrix'' approach in the spirit of Tutte's
decomposition.  An exact discrete expression for the two-point
function was given in \GEOD\ for planar quadrangulations and more
generally bipartite planar maps with a control on face degrees, whose
scaling form agrees with \AW. This exact expression makes use of a
coding of such maps by trees, leading to discrete recurrence equations
whose solution was guessed.

The purpose of this paper is to present an unexpected yet remarkable
connection between the two above notions showing that the information
about the distance in maps is actually hidden in the global problem of
enumerating maps with a boundary. More precisely, the moments and the
distance-dependent two-point function constitute two possible
expansions of the same quantity, the resolvent: the moments form its
power series expansion while the two-point function encodes its {\it
continued fraction expansion}.  Using the standard theory of continued
fractions this allows in particular to obtain explicit expressions for
the two-point function via the known techniques for computing the
moments. This program is carried out in detail in this paper. We begin by
an overview of our main results in Sect.~1.2, where precise
definitions for the moments and two-point function are given. The
actual connection between them is established in Sect.~2 by use of
the coding of maps by appropriate trees called mobiles. 
In Sect.~3 we derive an explicit
expression for the moments in terms of weighted paths by various
techniques: a direct non-constructive check, the solution of
Tutte's equation and a bijective construction. This
third approach also gives rise to a one-parameter family of expressions for
the moments forming so-called {\it conserved quantities}. Our
expression for moments is used in Sect.~4 to derive explicit
expressions for the two-point function in terms of Schur
functions. The special case of bipartite maps, which was the scope of
\GEOD, is discussed in Sect.~5, while the particularly simple cases
of triangulations and quadrangulations are addressed in Sect.~6.
Concluding remarks and discussion are gathered in Sect.~7.

\subsec{Overview of the main results}

A {\it planar map} is a connected graph (possibly with loops and
multiple edges) drawn on the sphere without edge crossings, and
considered up to continuous deformation. It is made of {\it vertices,
edges} and {\it faces}. The {\it degree} of a vertex or face is the
number of edges incident to it (counted with multiplicity). A map is
{\it rooted} if one of its edges is distinguished and oriented, the
{\it root face} being the face on the right of the root edge and the
{\it root degree} being the degree of the root face.

In this paper, we generally consider enumeration problems for rooted
planar maps subject to a control on face degrees, i.e. for each
positive integer $k$ we prescribe the number of faces with degree
$k$. In the language of generating functions, this is equivalent to
attaching a weight $g_k$ to each face of degree $k$, where
$(g_k)_{k\geq 1}$ is a sequence of indeterminates, so that the global
weight of a map is the product of the weights of its faces.
Setting $g_k=0$ for all odd
$k$ amounts to considering the planar maps which only have faces of
even degree, which are the {\it bipartite} planar maps. Drastic
simplifications occur in the bipartite case and these will be
discussed in due course. 

As a general preliminary remark, most quantities introduced in this
paper will be formal power series in the $g_k$'s (with integer or
rational coefficients) but, for the sake of concision, this will not
be apparent on the notations -- we shall write $X$ rather than
$X(g_1,g_2,g_3,\ldots)$ for instance. If extra variables are involved,
they shall be written explicitly.

The fundamental observation of this paper is a combinatorial identity
between some {\it a priori} unrelated families of generating
functions for rooted planar maps with the above weights. 
On the one hand, for each positive
integer $n$, we consider the generating function $F_n$ for rooted
planar maps with root degree $n$, which we call the $n$th {\it moment}. 
By convention we do not attach a
weight $g_n$ to the root face and we set $F_0=1$.  We may combine all
values of $n$ into the generating function $F(z)=\sum_{n=0}^{\infty}
F_n z^n$ where $z$ is an extra variable. This quantity essentially
coincides with the planar ``resolvent'' encountered in the context of
matrix integrals. On the other hand, again for
each positive integer $n$, we consider the generating function for
planar maps with two marked vertices at distance $n$, which we call
the distance-dependent {\it two-point function}. 
To properly root these planar maps, we take
as root an edge incident to and pointing away from one of the marked vertices.
Then the endpoint of the root edge must be at distance $n-1$, $n$ or
$n+1$ from the other marked vertex. We distinguish between these three
cases and denote by $r_n$, $t_n$ and $r_{n+1}$ the respective
generating function (clearly reversing the orientation of the root
edge allows to identify the third case with the first one for $n \to
n+1$). We furthermore introduce the ``cumulative'' generating
functions $R_n = \sum_{i=1}^n r_i$ and $T_n = \sum_{i=0}^n t_i$,
corresponding to two marked vertices at distance less than or equal
to $n$. For convenience, a conventional term $1$ is also included 
in $r_1$. Then, our observation is that
the sequences $(R_n)_{n \geq 1}$
and $(S_n)_{n \geq 0}=(\sqrt{T_n})_{n \geq 0}$ form the Jacobi type 
continued fraction (or {\it J-fraction}) expansion of $F(z)$, namely
\eqn\continued{F(z) = \sum_{n=0}^{\infty} F_n z^n =
{1 \over {\displaystyle 1 - S_0 z - {R_1 z^2
\over \displaystyle 1 - S_1 z - {R_2 z^2 \over 1 - \cdots}}}}.}
As discussed in detail in Sect.~2, this follows from the
correspondence between planar maps and some labeled trees called
mobiles \MOB, and the combinatorial theory of continued fractions
[\xref\FLAJFRAC,\xref\FLAJFRACBIS]. 
Calling $F_n$ the $n$th moment is consistent with the usual 
terminology of J-fractions \Wall. A classical result states that the sequence
$(R_n)_{n \geq 1}$ is closely related to the {\it Hankel determinants} of
moments, namely
\eqn\Rndet{R_n
= {H_{n} H_{n-2} \over H_{n-1}^2}, \qquad
  H_n = \det_{0 \leq i,j \leq n} F_{i+j}}
(with the convention $H_{-1}=1$) and similarly
the sequence $(S_n)_{n \geq 0}$ is given by {\it Hankel minors} through
\eqn\Sndet{\sum_{i=0}^{n} S_i = {\tilde{H}_n \over H_n}, \qquad
  \tilde{H}_n = \det_{0 \leq i,j \leq n} F_{i+j+\delta_{j,n}}.}

In the specific context of planar maps, we derive in Sect.~3 a general 
formula for the moments $F_n$ in terms of the weights $(g_k)_{k \geq 1}$.
A very peculiar structure emerges when this formula is
substituted into the determinants $H_n$ and $\tilde{H}_n$: after some
simple manipulations we recognize instances of a classical identity
for symplectic Schur functions, see Sect.~4 for details. Let us
now state the important results of Sects.~3 and 4.

Our expression for the moments involves generating functions for {\it
three-step paths}: these are lattice paths in the discrete Cartesian
plane which consist of {\it up-steps} $(1,1)$, {\it level-steps}
$(1,0)$ and {\it down-steps} $(1,-1)$.  We denote by $P(n;R,S)$ the
generating function for three-step paths going from $(0,0)$ to
$(n,0)$, where a weight $S$ is attached to each level-step, a weight
$\sqrt{R}$ is attached to each up- or down-step and the global weight
of a path is the product of the weights of its steps (clearly, this
global weight involves an integer power of $R$). Furthermore, we say
that a path is {\it positive} if it only visits vertices with
non-negative ordinate, and we denote by $P^+(n;R,S)$ the generating
function for positive three-steps paths from $(0,0)$ to $(n,0)$, also
known as {\it Motzkin paths} of length $n$, with the same weights.
We have the explicit expressions $P(n;R,S)= \sum_{i=0}^{\lfloor n/2 \rfloor}
{n! \over (i!)^2 (n-2i)!} R^i S^{n-2i}$ and $P^+(n;R,S)=
\sum_{i=0}^{\lfloor n/2 \rfloor} {n! \over i! (i+1)! (n-2i)!} R^i S^{n-2i}$.
Then, the general formula for the moments is  
\eqn\Fpaths{F_n = \sum_{q=0}^{\infty} A_q P^+(n+q;R,S)}
where the coefficients $A_q$ are given by
\eqn\apaths{A_q = R \left( \delta_{q,0} -
  \sum_{k=q+2}^{\infty} g_k P(k-q-2;R,S) \right)}
and where the step weights $R$, $S$ are themselves power series in
$(g_k)_{k \geq 1}$,
which are implicitly determined by
\eqn\RSeqs{S = \sum_{k=1}^{\infty} g_k P(k-1;R,S), \qquad
  R = 1 + {1 \over 2} \sum_{k=1}^{\infty} g_k P(k;R,S) - {S^2 \over 2}.} 

Now, substituting the expression \Fpaths\ into the Hankel determinant $H_n$,
elementary row and column manipulations, which amount to a natural
decomposition of lattice paths, lead to
\eqn\Hdetb{H_n = R^{n(n+1) \over 2} \det_{0 \leq i,j \leq n}
(B_{i-j} - B_{i+j+2}), \qquad B_i = \sum_{q=0}^{\infty} A_q P_i(q;R,S)}
where $P_i(q;R,S)$ is the generating function for three-step paths from
$(0,0)$ to $(q,i)$. Note that $P_{-i}(q;R,S)=P_i(q;R,S)$, hence 
$B_{-i}=B_i$. 
Similarly for $\tilde{H_n}$ we have
\eqn\tHdetb{\tilde{H}_n = (n+1) S H_n + R^{n^2 + n + 1 \over 2}
  \det_{0 \leq i,j \leq n} (B_{i-j-\delta_{j,n}} - B_{i+j+\delta_{j,n}+2}).} 
At this stage, we fix a positive integer $p$ and consider maps whose
faces have degree at most $p+2$, i.e. we set $g_k=0$ for $k>p+2$. This
implies, by the expression \apaths, that $A_q$ vanishes for $q>p$ and moreover,
since $P_i(q;R,S)=0$ for $|i|>q$, that $B_i$ vanishes for $i>p$. The
matrices appearing in \Hdetb\ and \tHdetb\ are therefore band matrices.
Their determinants can be identified, up to a normalization, with 
characters of the symplectic
group ${\rm Sp_{2p}}$ \FULHA, also known as {\it symplectic Schur
functions}. Following \FULKRA, we denote by ${\rm
sp}_{2p}(\lambda,{\bf x})$ the symplectic Schur function associated with
the partition $\lambda$, which is a symmetric function of $2p$
variables ${\bf x}=(x_1,1/x_1,x_2,1/x_2,\ldots,x_p,1/x_p)$. Here, these
variables are the $2p$ roots of the ``characteristic equation''
\eqn\chareqb{\sum_{n=-p}^{p} B_{n} x^n = 0}
(recall that $B_{-n}=B_n$, hence the equation is invariant under $x\to 1/x$).
By the expressions \apaths\ for $A_q$, \Hdetb\ for $B_i$
and the identity 
$\sum_{i=-q}^q P_i(q;R,S)\, x^i= (\sqrt{R}\, x+ S+\sqrt{R}/x)^q$,
the characteristic equation reads explicitly
\eqn\chareqP{1=\sum_{k=2}^{p+2} g_k \sum_{q=0}^{k-2} P(k-2-q;R,S)
\left(\sqrt{R}\, x+S+{\sqrt{R}\over x}\right)^q\ .}
The partitions associated with the 
determinants appearing in \Hdetb\ and \tHdetb\ 
are respectively the
``rectangular'' partition $\lambda_{p,n+1}=(n+1)^p$ made of $p$ parts of size
$n+1$, and the ``nearly-rectangular'' partition
$\tilde{\lambda}_{p,n+1}=(n+1)^{p-1}n$ 
made of $p-1$ parts of size $n+1$ plus one
part of size $n$. In summary, we find that
$H_n\propto {\rm sp}_{2p}(\lambda_{p,n+1},{\bf x})$ and 
${\tilde H}_n-(n+1)S H_n \propto 
{\rm sp}_{2p}(\tilde{\lambda}_{p,n+1},{\bf x})$. 
Using the Weyl character formula for ${\rm sp}_{2p}(\lambda,{\bf x})$, we 
end up with the following ``nice'' formulas for $R_n$ and $S_n$, which involve 
determinants of size $p$, independently of $n$:
\eqn\finalRn{R_n=R\  {
\det\limits_{1 \leq i,j \leq p} (x_i^{n+1+j} - x_i^{-n-1-j}) 
 \det\limits_{1 \leq i,j \leq p} (x_i^{n-1+j} - x_i^{-n+1-j}) 
\over \left( \det\limits_{1 \leq i,j \leq p} (x_i^{n+j} -
x_i^{-n-j}) \right)^2 }\ ,}
\eqn\Snfinal{S_n=S-\sqrt{R} \left( 
{\det\limits_{1 \leq i,j \leq p} (x_i^{n+1+j-\delta_{j,1}} -
x_i^{-n-1-j+\delta_{j,1}}) \over
\det\limits_{1 \leq i,j \leq p} (x_i^{n+1+j} - x_i^{-n-1-j})}
- 
{\det\limits_{1 \leq i,j \leq p} (x_i^{n+j-\delta_{j,1}} -
x_i^{-n-j+\delta_{j,1}}) \over
\det\limits_{1 \leq i,j \leq p} (x_i^{n+j} -
x_i^{-n-j})} \right)\ .}

A number of simplifications occur in the bipartite case, discussed
in detail in Sect.~5. Clearly $F_n=0$ for 
$n$ odd and the continued fraction expansion of $F(z)$ naturally reduces
to the Stieltjes type 
\eqn\Stieltjes{F(z) = \sum_{n=0}^{\infty} F_{2n} z^{2n} =
{1 \over {\displaystyle 1 - {R_1 z^2
\over \displaystyle 1 - {R_2 z^2 \over 1 - \cdots}}}}.}
Consistently, the quantities $S$, $\tilde{H}_n$ and $S_n$ vanish while 
the Hankel determinants factorize as 
\eqn\factorhank{H_n= h^{(0)}_{\left\lfloor {n\over 2}\right\rfloor}
h^{(1)}_{\left\lfloor {n-1\over 2}\right\rfloor}, \qquad
h^{(0)}_n=\det_{0\leq i,j\leq n} F_{2i+2j}, \quad
h^{(1)}_n=\det_{0\leq i,j\leq n} F_{2i+2j+2}\ .
}
This factorization property is also apparent at the level 
of Schur functions. Indeed, for bipartite maps whose faces have degree
at most $2p+2$, we have a characteristic equation of the form
\eqn\charbip{\sum_{n=-p}^p B_{2n}x^{2n}=0} 
whose roots are of the form $(x_1,1/x_1,-x_1,-1/x_1, \ldots, x_p,1/x_p, -x_p,
-1/x_p)$, and ${\rm sp}_{4p}(\lambda,{\bf x})$ can be factorized as
the product of two symmetric functions of the variables $(x_1^2, 1/x_1^2,\ldots
x_p^2,1/x_p^2)$ (a symplectic and an odd orthogonal Schur function). In the end,
$R_n$ can be expressed in terms of determinants of size $p$ rather than
$2p$, which matches the expression found in \GEOD. 

The cases of triangulations and quadrangulations are particularly simple 
as the roots of their characteristic equation are of the respective forms
$(x,1/x)$ and $(x,1/x,-x,-1/x)$. Consequently, the two-point functions
in these cases can be expressed in terms of a single quantity $x$ solution of
some algebraic equation. This is discussed in Sect.~6 and a simple 
interpretation via one-dimensional hard dimers is given.

\newsec{From planar maps to continued fractions, via mobiles}

The main purpose of this section is to establish that the generating
functions for planar maps $F_n$, $R_n$ and $S_n$ defined in Sect.~1
are related by the J-fraction expansion \continued. The plan is as
follows. We recall some general results of the combinatorial theory of
continued fractions (Sect.~2.1), before reviewing the coding of
planar maps by mobiles (Sect.~2.2). We then interpret $F_n$, $R_n$
and $S_n$ as generating functions for mobiles and complete our proof
(Sect.~2.3). Finally we mention a few other outcomes of our approach
(Sect.~2.4).

\subsec{Reminders of the combinatorial theory of continued fractions}

In his seminal paper \FLAJFRAC, Flajolet gave a combinatorial
interpretation of continued fractions in terms of Motzkin paths. Let
us recall a few important results of this theory.  For the purposes of
this subsection, $(R_m)_{m\geq 1}$ and $(S_m)_{m\geq 0}$ may denote
arbitrary sequences of elements of a commutative ring. That is to say, 
we may temporarily forget about their map-related definition given in
Sect.~1.  We consider Motzkin paths, and more generally
positive three-step paths as defined above, with the following {\it
ordinate-dependent} weights: for all $m$,
\item{-} each down-step of the form $(t,m) \to (t+1,m-1)$ receives a weight $R_m$ ($m \geq 1$)
\item{-} each level-step of the form $(t,m) \to (t+1,m)$ receives a weight $S_m$ ($m \geq 0$)
\item{-} each up-step receives a weight $1$.
\par
\noindent Note that attaching ordinate-dependent weights to up-steps
would add essentially no generality. Let us now define $F_n$ as the
generating function for Motzkin paths of length $n$ with
ordinate-dependent weights (we also temporarily forget about the
map-related definition of $F_n$) and, attaching a further weight $z$
per step (of any type), $F(z)=\sum_{n=0}^{\infty} F_n z^n$ the
generating function for Motzkin paths of arbitrary length. Then, the
Continued Fraction Theorem (as named in \FLAJFRACBIS) states that
$F(z)$ is given by the J-fraction \continued.

In order to recall of its derivation, let us introduce a few notations which
also prove useful later. Given arbitrary non-negative integers $d$,
$d'$ and $n$, we denote by $Z_{d,d'}(n)$ the generating
function for positive three-step paths from $(0,d)$ to
$(n,d')$, with the above ordinate-dependent weights.
Furthermore, we denote by $Z^+_{d,d'}(n)$ the generating
function for such paths restricted to have all their ordinates larger
than or equal to $\min(d,d')$ (instead of 0). Clearly
$F_n=Z_{0,0}(n)=Z^+_{0,0}(n)$. Now, for the case $d=d'$ we have the relation
\eqn\decomppath{\sum_{n\geq 0} Z^+_{d,d}(n)\, z^n= {1\over
 1-S_{d}\, z-R_{d+1}\, z^2 \sum\limits_{n\geq 0} Z^+_{d+1,d+1}(n)\, z^n}}
which simply translates the ``arch decomposition'', namely by cutting
a path contributing to the left-hand side at each occurrence of the
ordinate $d$, it is bijectively decomposed into a sequence of two
types of objects: either (i) level-steps at ordinate $d$, weighted by
$S_d$, or (ii) arches made of the concatenation of an up-step, a
restricted path starting and ending at ordinate $d+1$ and a down-step,
with an overall weight $R_{d+1}Z_{d+1,d+1}(n')$ for some $n'$. By iterating
\decomppath\ starting at $d=0$, we deduce \continued. Note that the
J-fraction is a well-defined power series in $z$, because computing
$F_n$ only requires $\lfloor n/2 \rfloor$ iterations (as Motzkin paths
of length $n$ reach at most the ordinate $\lfloor n/2 \rfloor$). If we
iterate \decomppath\ starting at an arbitrary $d$, we obtain
\eqn\truncation{\sum_{n\geq 0} Z^+_{d,d}(n)\, z^n =
{1 \over {\displaystyle 1 - S_d z - {R_{d+1} z^2
\over \displaystyle 1 - S_{d+1} z - {R_{d+2} z^2 \over 1 - \cdots}}}}}
where the right-hand side is called a {\it truncation} of the fundamental
fraction \continued. Many other identities are known 
(e.g for $\sum_{n\geq 0} Z^+_{d,d'}(n)\, z^n$ and $\sum_{n\geq 0}
Z_{d,d'}(n)\, z^n$ via ``last-passages decompositions'') but we shall
not need them here.

A particular case is when the weights do not depend on the ordinate,
i.e. we set $R_m=R$ for all $m \geq 1$ and $S_m=S$ for all $m \geq
0$. Then, we recover the three-step path generating functions defined
in Sect.~1, namely $Z_{d,d'}(n)= R^{(d-d')/2} P_{d'-d}(n;R,S)$ and
$Z^+_{d,d'}(n)= R^{(d-d')/2} P^+_{d'-d}(n;R,S)$. Note that the 
$R^{(d-d')/2}$ prefactor arises from the slightly different weighting
convention. Most notably, we have
\eqn\motzfrac{\sum_{n\geq 0} P^+(n;R,S)\, z^n = {1 \over
{\displaystyle 1 - S z - {R z^2 \over \displaystyle 1 - S z - {R z^2
\over 1 - \cdots}}}}.}
Clearly this series is the solution of the quadratic equation $X = 1 +
S z X + R z^2 X^2$. More generally, if we only assume that $R_m \to R$
and $S_m \to S$ for $m \to \infty$ (in any sense of convergence),
then obviously $Z_{d,d}(n) \to P(n;R,S)$ and
$Z^+_{d,d}(n) \to P^+(n;R,S)$ for $d \to \infty$.

Let us now consider the inverse problem of determining the sequences
$(R_m)_{m\geq 1}$ and $(S_m)_{m\geq 0}$ in the J-fraction expansion,
knowing the power series expansion of $F(z)$ i.e. the sequence
$(F_n)_{n \geq 0}$. It is a classical result
that this problem is solved using Hankel determinants, and the
solution is given by relations \Rndet\ and \Sndet. See for
instance the beautiful proof by Viennot \Viennot\ based on a combinatorial 
interpretation involving configurations of non-intersecting Motzkin paths.

\subsec{Review of the coding of planar maps by mobiles}

\fig{The mobile construction applied (a) to an edge of type $(n,n-1)$
and (b) to an edge of type $(n,n)$. In (a) a mobile edge (thick line)
connects an unlabeled vertex (big dot) to a labeled one. In (b) a
bivalent flagged vertex (lozenge) is created and connected to both
adjacent unlabeled vertices.}{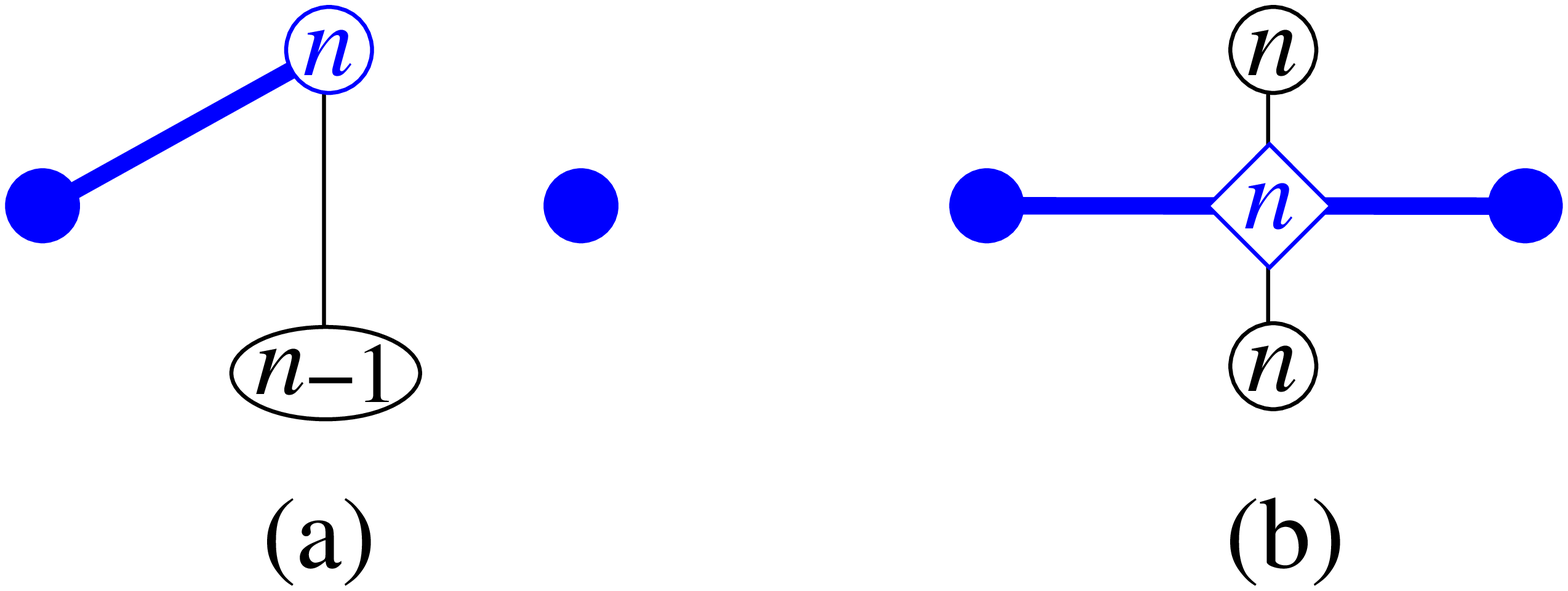}{10.cm}
\figlabel\edgerule

We now return to the context of planar maps and review their coding by
mobiles presented in Ref.~\MOB. The coding applies to {\it pointed}
maps, i.e. maps with a distinguished vertex called the origin, and
works as follows: we start by labeling each vertex of the map by a
non-negative integer equal to its graph distance from the origin and
add a new unlabeled vertex at the center of each face. With the above
labeling, any edge of the map is either of type $(n,n-1)$, i.e.
connects a vertex labeled $n$ to a vertex labeled $n-1$ for some
$n\geq 1$, or of type $(n,n)$ for some $n\geq 0$. For each edge of
type $(n,n-1)$, we draw a new edge connecting its incident vertex
labeled $n$ to the unlabeled vertex sitting at the center of the face
on the right of the edge (oriented from $n$ to $n-1$). This procedure
is displayed in Fig.~\edgerule-(a). For each edge of type $(n,n)$, we
add in its middle a ``flagged'' vertex with flag $n$ and connect it by
two new edges to the two unlabeled vertices at the center of the two
incident faces, as displayed in Fig.~\edgerule-(b). It was shown in
\MOB\ that the new edges form a tree which spans all the labeled
vertices with label $n\geq 1$ (i.e. all the original vertices of the
map but its origin), as well as all the added unlabeled vertices and
all the added flagged vertices.

\fig{An illustration of the property (P) around an unlabeled vertex.}{
  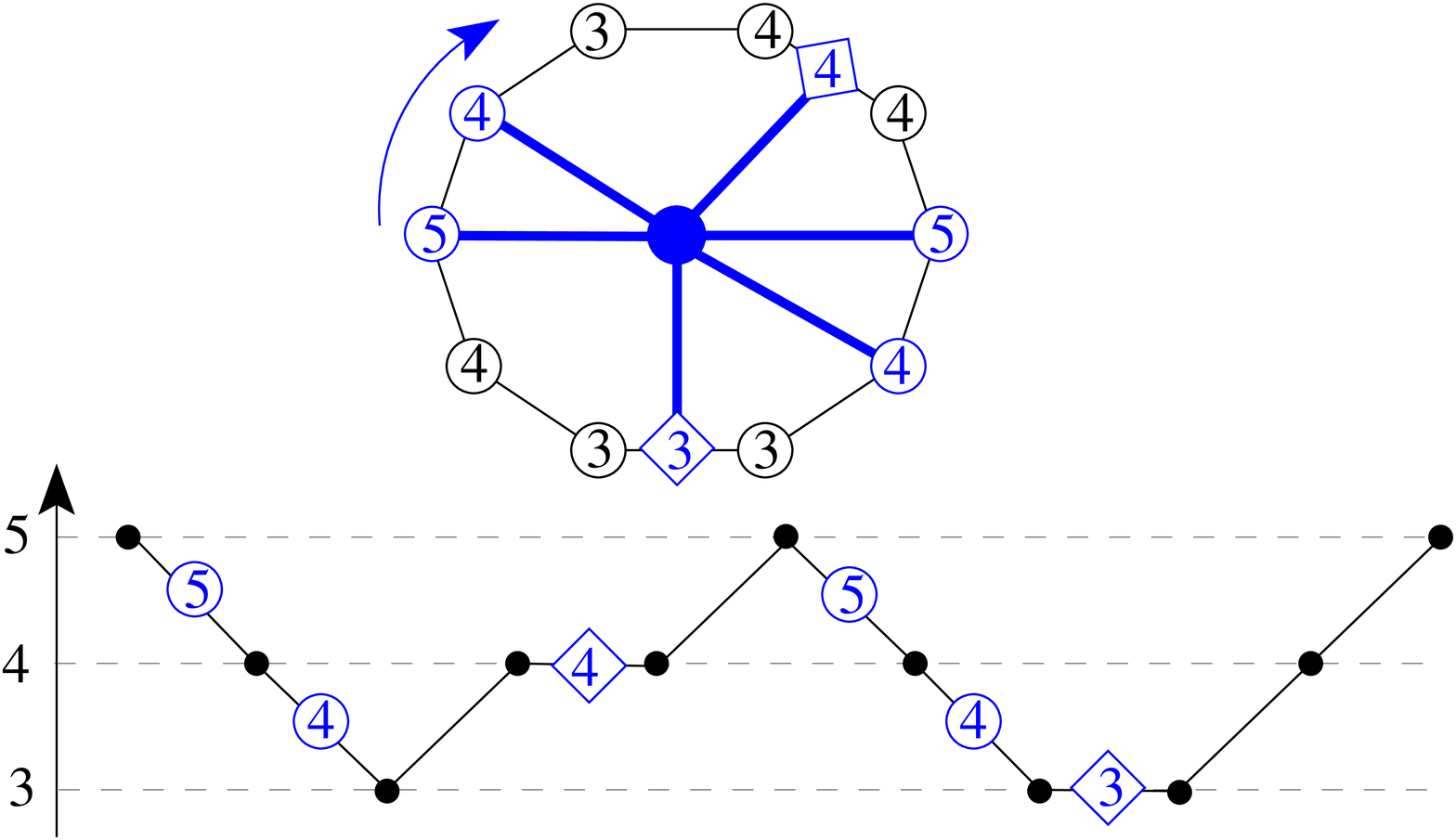}{10.cm}
\figlabel\facerule

The labels and flags satisfy the property:
\item{(P)} for each unlabeled vertex, the {\it clockwise} sequence
of labels and flags on its adjacent vertices matches exactly the sequence
of labels and flags obtained by attaching a label $n$ to each down-step 
$(t,n)\to (t+1,n-1)$ and a flag $n$ to each level-step $(t,n)\to (t+1,n)$ 
of some three-step path with identical initial and final ordinates.
\par\noindent Note that the three-step path in (P) is unique up to
cyclic shifts, and that its length is given by the number of flags
plus twice the number of labels. Returning to the original map, the
sequence of ordinates in the three-step path simply reproduces the
clockwise sequence of distances from the origin to the vertices around
the face associated with the unlabeled vertex at hand, while the
length of the path is nothing but the degree of this face (see
Fig.~\facerule).

In all generality, a {\it mobile} is defined as a plane tree with
three types of vertices: unlabeled vertices, labeled vertices carrying
an arbitrary integer label and flagged vertices carrying an arbitrary
integer flag, with edges connecting only unlabeled vertices to labeled
or flagged ones in such a way that the flagged vertices have degree
$2$ and property (P) holds around each unlabeled vertex. Note that
shifting all labels and flags in a mobile by a fixed integer preserves
property (P) so that the resulting object is still a mobile.

The above coding associates to each pointed map a mobile. As shown in
\MOB, it is a bijection between pointed planar maps and mobiles 
satisfying the two extra requirements:
\item{(R1)} labels and flags are respectively positive and non-negative;
\item{(R2)} there is at least one label $1$ or one flag $0$.
\par\noindent Strictly speaking, the tree reduced to a single
vertex labeled $1$ satisfies (R1) and (R2), and corresponds to a
conventional {\it empty map}.

Furthermore, the coding yields the following one-to-one
correspondences between elements of a map and of its associated
mobile.  Map vertices at distance $n\geq 1$ from the origin in the map
correspond to mobile vertices labeled $n$. Map edges of type $(n,n)$,
i.e. connecting two vertices at distance $n$ from the origin,
correspond to flagged vertices with flag $n$ ($n\geq 0$) while map
edges of type $(n,n-1)$ correspond to mobile edges incident to a
labeled vertex with label $n$ ($n\geq 1$). Faces of the map correspond
to unlabeled vertices of the mobile and moreover, the clockwise
sequence of distances from the origin of the vertices incident to a
given face is directly read off the ordinates of the (unique cyclic)
three-step path ensuring property (P) around the corresponding
unlabeled vertex.

The mobile coding a pointed map is unrooted, i.e. it has no
distinguished edge nor vertex. We define conventionally a {\it rooted
mobile} as a mobile with a distinguished edge incident to a labeled
vertex (the other incident vertex being necessarily unlabeled), whose
label is then called the {\it root label}. The mobile reduced to a
single labeled vertex is considered as a rooted mobile even if it has
no edge. Regarding flagged vertices, we find more convenient to
introduce the notion of {\it half-mobile}, which is defined exactly as
a mobile except that it has one particular flagged vertex of degree
$1$ (the other flagged vertices being of degree $2$ as before) whose
flag is called the {\it root flag}.  A mobile having a distinguished
edge incident to a flagged vertex with flag $n$ can clearly be seen as
a pair of half-mobiles having the same root flag $n$.  Rooted mobiles
and pairs of half-mobiles code naturally pointed rooted maps, namely
maps with both a marked vertex (the origin) and a marked oriented edge
(the root edge). More precisely, a pointed rooted map of type $(n \to
n-1)$, i.e. whose root edge points from a vertex at distance $n$ from
the origin to a vertex at distance $n-1$, is coded bijectively by a
rooted mobile satisfying (R1)-(R2) and having root label $n$ (we
distinguish the mobile edge associated with the root edge), while a
pointed rooted map of type $(n \to n)$ is coded bijectively by a pair
of half-mobiles with root flag $n$ and whose union also satisfies
(R1)-(R2) (we split the mobile at the flagged vertex associated with
the root edge, whose orientation allows to distinguish the two
half-mobiles). For the last remaining type $(n-1 \to n)$ we may simply
reverse the root edge orientation. This encompasses all possible
types of the root edge, when $n$ varies. Note that rooted maps can be
seen as pointed rooted maps of type $(1 \to 0)$ or $(0 \to 0)$.

It is worth noting that a bijection also exists with mobiles where the
conditions (R1) and (R2) are waived. It is obtained by now considering
pointed rooted maps where the distance between the origin of the map
and the starting point of the root edge is no longer fixed. 
Denoting by $d$ the distance of the extremity of
the root edge farthest from the origin, we now label each vertex by
its distance from the origin minus $d$ so that the pointed rooted map
is either of type $(0 \to 0)$, $(0 \to -1)$ or $(-1 \to 0)$. Applying
the mobile construction rules of Fig.~\edgerule\ to all edges of the
map creates a mobile where the condition (P) still holds but where the
conditions (R1) and (R2) have been waived. For type $(0 \to 0)$ we
finally obtain a pair of half-mobiles with root flag $0$, while for
type $(0 \to -1)$ or $(-1 \to 0)$ we obtain a rooted mobile with root
label $0$ (for a proper bijection we need to adjoin a sign $\pm 1$ to
the rooted mobile, in order to keep track of the orientation of the
root edge).

\subsec{Mobile interpretation of the generating functions $F_n$, $R_n$, $S_n$}

We now want to show that \continued\ holds, when $F_n$, $R_n$ and
$S_n$ are the generating functions for planar maps defined in Sect.~1.
This requires identifying these quantities as generating functions
for mobiles, and $F_n$ will naturally appear as a sum over Motzkin
paths with weights $R_m$ per down-step starting at ordinate $m$ and
$S_m$ per level-step at ordinate $m$. As a general preliminary remark,
note that the weight $g_k$ per face of degree $k$ in the map
translates, in the mobile language, into a weight $g_k$ per unlabeled
vertex whose three-step path associated via (P) has length $k$.

To begin with, let us discuss the case of $R_n$ ($n \geq 1$). From our
definition of Sect.~1.2 and the above discussion, $R_n=\sum_{i=1}^n
r_i$ enumerates pointed rooted maps of type $(i \to i-1)$ for some
$i$ between $1$ and $n$, which are in bijection with rooted mobiles
satisfying (R1)-(R2) and having a root label $i$ between $1$ and
$n$. These mobiles are in turn in bijection with rooted mobiles {\it
satisfying} (R1) {\it only} and having root label $n$. Indeed, we may
transform the label $i$ in the first family of mobiles into a label
$n$ by shifting all labels and flags of the mobile by $n-i$. The
condition (R1) is preserved by this non-negative shift but (R2) is no
longer valid if $i<n$. Conversely, assuming that (R1) is satisfied in
a rooted mobile with root label $n\geq 1$, the minimal label among
labeled vertices necessarily lies between $1$ and $n$, i.e. is of the
form $n-i+1$ for some $i$ between $1$ and $n$.  From the general
condition (P) around each unlabeled vertex, we easily deduce that the
minimal flag in the mobile is at least $n-i$. Shifting all labels and
flags by the (non-positive) quantity $i-n$ preserves (R1) and restores
(R2), while it transforms the root label $n$ into an $i$ between $1$
and $n$.  To summarize, $R_n$ may be identified as the generating
function for rooted mobiles satisfying (R1) only and having root label
$n$. In this identification, the conventional term $1$ added in $r_1$,
hence in $R_n$, accounts for the mobile reduced to a single labeled
vertex with label $n$.

As for $T_n=\sum_{i=0}^n t_i$, it enumerates pointed rooted maps of
type $(i \to i)$ for some $i$ between $0$ and $n$, which are in bijection
with pairs of half-mobiles with root flag $i$ and whose union
satisfies (R1)-(R2). By a similar argument as above, these are in
bijection with pairs of half-mobiles with root flag $n$ satisfying
(R1) only. Note that (R1) independently applies on each
half-mobile. We arrive at $T_n=(S_n)^2$ where we identify directly
$S_n$ as the generating function for half-mobiles satisfying (R1) only
and with root flag $n$.

\fig{The generating function $F_n$ for rooted maps with root degree $n$
is also that $Z_{0,0}(n)$ for Motzkin paths of length $n$ with
ordinate-dependent weights.}{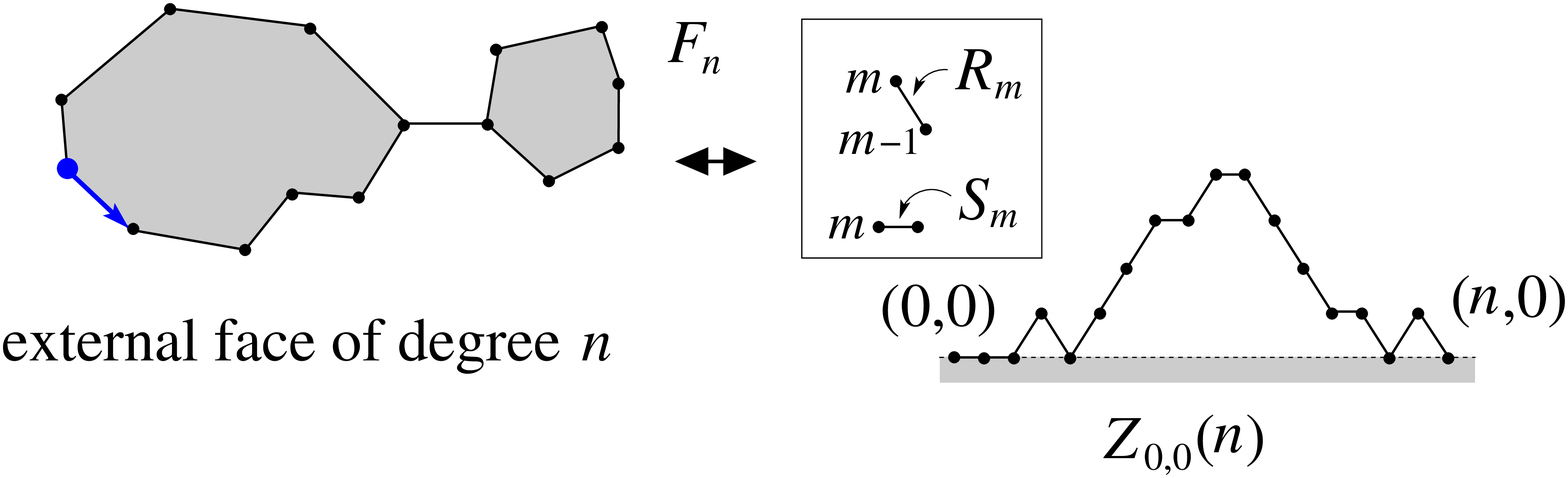}{13.cm}
\figlabel\motzmapzero

Let us now come to the mobile interpretation of $F_n$. From its
definition of Sect.~1.2, it enumerates rooted maps with a root face
of degree $n$.  Taking as origin of the map the origin of the root
edge, we may look at the clockwise sequence of distances from this
origin of the successive vertices incident to the root face in the
map. For convenience, in the planar representation, we choose as
external face the root face itself so that the clockwise orientation
around it corresponds in practice to the counterclockwise orientation
around the rest of the map. Starting from the origin, this sequence of
distances forms a Motzkin path of length $n$.  The down-steps and the
level-steps of the Moztkin path correspond respectively to the labeled
vertices and the flagged vertices connected to the unlabeled vertex
associated with the root face in the mobile. By removing this
unlabeled vertex and all its incident edges, we obtain a sequence of
rooted mobiles and half-mobiles (when removing an edge incident to a
labeled vertex, we keep track of its position by distinguishing, say,
the next edge incident to the same labeled vertex in clockwise
direction). More precisely, on the Motzkin path, a down-step of the
form $(t,m)\to (t+1,m-1)$ corresponds to a rooted mobile with root
label $m$ satisfying (R1) (possibly reduced to a single labeled
vertex), while a level-step $(t,m)\to (t+1,m)$ corresponds to a
half-mobile with root flag $m$ also satisfying (R1).  Conversely,
starting from a Motzkin path and a sequence of rooted and half-
mobiles so associated with its down- and level-steps, a complete
mobile is immediately obtained by connecting the mobiles to a new
unlabeled vertex, and properties (P)-(R1)-(R2) are satisfied by
construction.  This clearly forms a bijection. Translated in the
language of generating functions, it says that $F_n$ is equal to the
generating function for Motzkin paths of length $n$ with
ordinate-dependent weights $(R_m)_{m \geq 1}$ and $(S_m)_{m \geq 0}$ (see
Fig.~\motzmapzero), i.e. $F_n = Z_{0,0}(n)$ with the notations of
Sect.~2.1. Equation \continued\ follows from the general
combinatorial theory of continued fractions.

To conclude this section, we mention that the continued fraction
expansion may be directly interpreted at the level of maps, with no
recourse to the coding by mobiles. It makes use of a particular
decomposition of pointed rooted maps into slices, as explained in
Appendix A.

\subsec{Related results}

Before proceeding to the next section, let us list a few other
enumerative consequences of the coding of maps by mobiles. First,
$R_n$ and $S_n$ satisfy recursive equations which translate
decompositions of the corresponding mobiles. Indeed, let us consider a
rooted (resp.\ half-) mobile with root label (resp.\ flag) $n$
satisfying (R1), not reduced to a single labeled vertex. We now
decompose the mobile around the unlabeled vertex incident to the
distinguished edge (resp.\ adjacent to the root flagged vertex). The
three-step path associated with this unlabeled vertex has an arbitrary
length $k \geq 1$ (thus leads to a weight $g_k$), is positive but does
not necessarily attain 0. The distinguished edge yields a
distinguished down-step (resp.\ level-step) starting from ordinate
$n$. Hence all the other steps form a three-step path of length $k-1$
starting at ordinate $n-1$ (resp.\ $n$) and ending at ordinate $n$. As
in the above discussion of $F_n$, each down- or level-step starting at
an ordinate $m$ is associated with a rooted or half-mobile with root
label or flag $m$ satisfying (R1) (except for the distinguished
level-step when decomposing a half-mobile), and the decomposition is
bijective. In the end, this translates into the relations
\eqn\rnsnrecur{\left\{ \eqalign{
  R_n &= 1 + \sum_{k=1}^{\infty} g_k R_n Z_{n-1,n}(k-1) =
  1 + \sum_{k=1}^{\infty} g_k Z_{n,n-1}(k-1) \cr
  S_n &= \sum_{k=1}^{\infty} g_k Z_{n,n}(k-1)} \right.}
where $1$ stands for the mobile reduced to a single vertex labeled
$n$, and the $Z$'s are defined as in Sect.~2.1. In the second equation 
for $R_n$, we used the identity $Z_{n,n-1}(k-1)=R_n Z_{n-1,n}(k-1)$ obtained
by reflecting the paths with respect to a vertical axis. Note that the $Z$'s
depend implicitly on the sequences $(R_m)_{m \geq 1}$ and $(S_m)_{m
\geq 0}$ hence relations \rnsnrecur\ are of recursive nature. It is easily
seen that they uniquely determine all $R_n$'s and $S_n$'s as power
series in the $g_k$'s. These equations were already derived in \MOB\
and also in \PDFRaman\ using a different coding of maps (the bipartite
case was first discussed in \GEOD\ where an explicit solution was guessed).

Let us now briefly discuss the simpler case of mobiles in which
condition (R1) is waived. By a simple shift of labels by $-n$, $R_n$
(resp.\ $S_n$) enumerates rooted mobiles (resp.\ half-mobiles) with
root label (resp.\ flag) $0$ and with their labels all strictly larger
than $-n$ and their flags all larger than or equal to $-n$. Sending $n
\to \infty$ amounts to waiving this lower bound. This implies that
$R_n$ and $S_n$ converge for $n \to \infty$ in the sense of power
series, i.e. each of their coefficients stabilizes. We denote by $R$
and $S$ their respective limits, which are nothing but generating
functions for respectively rooted mobiles with root label $0$ and
half-mobiles with root flag $0$. By the discussion at the end of the
Sect.~2.2, we find that $2R+S^2$ is the generating function
for pointed rooted maps.  From \rnsnrecur\ we deduce
\eqn\rsrecurbis{\left\{ \eqalign{
  R &= 1 + \sum_{k=1}^{\infty} g_k \sqrt{R} \, P_{-1}(k-1;R,S) \cr
  S &= \sum_{k=1}^{\infty} g_k P(k-1;R,S)} \right.}
with $P_{-1}$ and $P$ defined as in Sect.~1. Note that these
relations are tantamount to \RSeqs\ by virtue of the identity
$P(k;R,S)=S P(k-1;R,S) + 2 \sqrt{R} P_{-1}(k-1;R,S)$. Note also that the
$k=1$ term in the first line of \rsrecurbis\ vanishes so that we may start the
summation at $k=2$ for $R$.

\fig{The generating function $F_n^{\bullet}$ for pointed rooted maps
with root degree $n$ is also that $P(n;R,S)$ for three-step paths of
length $n$ starting and ending at ordinate $0$, with
ordinate-independent weights.}{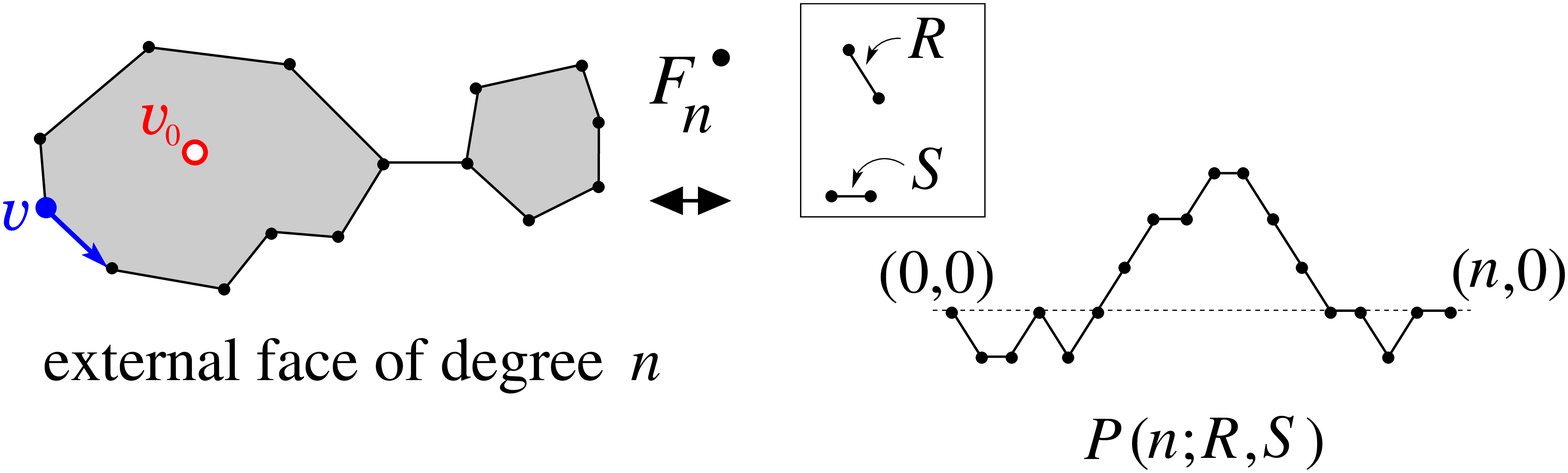}{13.cm}
\figlabel\motzmapfour

We finally consider pointed rooted maps whose root face has a
prescribed degree $n$, and denote by $F_n^{\bullet}$ the associated
generating function (where the root face does not receive a weight
$g_n$). $F_n^{\bullet}$ only differs from $F_n$ by the extra marking
of a vertex or, said otherwise, we no longer impose that the origin of
the map be the starting point of the root edge. Denoting by $d$ the
distance from the origin of the map to the origin of the root edge, we
now label each vertex by its distance from the origin of the map minus
$d$. Applying the mobile construction rules of Fig.~\edgerule, we
obtain a mobile with a distinguished unlabeled vertex corresponding to
the root face. Its associated three-step path codes the distances 
(minus $d$) from
the origin to the map vertices incident to the root face. The origin
of the root edge distinguishes a step starting at ordinate $0$ on the
path, therefore the path can be seen as a three-step path of length
$n$ which starts and ends at ordinate $0$, but which is not necessarily
positive. Decomposing around the distinguished unlabeled vertex, each
down- or level-step of the path yields a rooted or half-mobile, whose
root label or flag we may set to $0$ by a suitable shift. This
decomposition is bijective, and translates into the relation (see
Fig.~\motzmapfour)
\eqn\Fnpoint{F_n^{\bullet} = P(n;R,S).}
This relation is analogous to the relation $F_n=Z_{0,0}(n)$
established above, yet it is considerably simpler because it involves
three-step paths with weights independent from the ordinate. In this
sense, it is of the same nature as formula \Fpaths, which we shall now prove.

\newsec{An explicit expression for the moments}

We now turn to the derivation of formula \Fpaths\ for the moments
$F_n$ as defined in Sect.~1. In contrast with the J-fraction
expansion \continued, which follows from a general scheme where planar
maps and mobiles turn out to fit, formula \Fpaths\ seems deeply
related to the map structure. One of our purposes is to show that it
naturally results from several approaches. First, we discuss a direct
yet non-constructive check based on the comparison with formula
\Fnpoint\ (Sect.~3.1). Then, we explain how formula \Fpaths\ arises
from the solution of Tutte's equation, which is the original approach
to enumerating maps (Sect.~3.2). Finally we present a bijective
proof involving distance-dependent generating functions and mobiles
(Sect.~3.3). Beyond our main purpose, this third approach yields a
one-parameter family of expressions for $F_n$, the so-called
``conserved quantities'', interesting on their own.

\subsec{Proof via pointed rooted maps}

It has been noted on several occasions \CMS\ that considering maps
that are both pointed and rooted brings some simplifications. This is
illustrated by formula \Fnpoint\ in our context. A natural idea is to
deduce formula \Fpaths\ from it. Because $F_n^{\bullet}$ only differs
from $F_n$ by the extra marking of a vertex, there is a direct
relation between these two quantities.

To make this relation explicit, it is most convenient to introduce an
extra weight $u$ per vertex. We denote by $F_n^{\bullet}(u)$ and
$F_n(u)$ the correspondingly modified generating functions. (Note that
$u$ is actually a redundant parameter, as by a simple counting argument and
Euler's relation we have $F_n(u) = u^{n/2+1} \left. F_n(1)
\right|_{g_k \to g_k u^{k/2-1}}$ and the same for $F_n^{\bullet}$.)
On the one hand, relation \Fnpoint\ immediately generalizes as
\eqn\Fnpointu{F_n^{\bullet}(u) = u\, P(n;R(u),S(u)).}
where $R(u)$ and $S(u)$ enumerate mobiles and half-mobiles with an
additional weight $u$ per labeled vertex, and the extra factor $u$
accounts for the origin. On the other hand, the extra marking of a vertex
yields
\eqn\Fnudu{F_n^{\bullet}(u) = u {d \over du} F_n(u).}
Therefore $F_n(u)$ is a primitive of $P(n;R(u),S(u))$ with respect to the
variable $u$. It is shown in Appendix B that it has the explicit form
\eqn\Fnu{F_n(u) = \sum_{q=0}^{\infty} A_q(u) P^+(n+q;R(u),S(u))}
with $A_q(u)$ defined as in \apaths\ with $R,S$ replaced by $R(u),S(u)$.
By specializing to $u=1$ we recover the expression $\Fpaths$ for $F_n=F_n(1)$.

\subsec{Tutte's equation and its solution}

The generating function $F_n$ for rooted planar maps with root degree
$n$ is of utmost importance in Tutte's original approach \TUTEQ. Indeed,
the whole family $(F_n)_{n \geq 1}$ is uniquely determined as a power
series in the variables $(g_k)_{k \geq 0}$ by the equation
\eqn\tutteqn{F_n = \sum_{i=0}^{n-2} F_i F_{n-2-i} + \sum_{k=1}^{\infty}
  g_k F_{n+k-2}}
valid for all $n \geq 1$, with the convention $F_0=1$. This equation
directly expresses that, removing the root edge of a planar map with
root degree $n$, two situations may arise: either the map is split
into two connected components, which can be seen as two rooted planar
maps whose root degrees add up to $n-2$ (possibly one of these maps is
reduced to a single vertex and has root degree $0$), or the map is
not split, then its root degree is increased by $k-2$ where $k$ is the
degree of the face formerly on the left of the root edge. Passing to
$F(z)$ we get Tutte's equation \TUTEQ
\eqn\tutteqz{F(z) = 1 + z^2 F(z)^2 +
\sum_{k=1}^{\infty} g_k z^{2-k} \left( F(z) - \sum_{j=0}^{k-2}
z^j F_j \right)}
also known as {\it loop equation} in the context of matrix models.

The solution of Tutte's equation \tutteqz\ was obtained by Bender and
Canfield \BECA\ (with slight restrictions on the $g_k$), see also
\BMJ\ for a more general formulation and \CENSUS\ for the matrix model
counterpart. These authors were ultimately interested in the series
$F_2$ which, upon squeezing the bivalent root face, yields the
``true'' generating function for rooted planar maps. Our purpose here
is to show that the general formula \Fpaths\ for $F_n$ follows
straightforwardly from their discussion, which we now recall briefly.

We first
observe that Eq.~\tutteqz\ is quadratic in $F(z)$ (viewing the
terms $\sum_{j=0}^{k-2} z^j F_j$ as constants), hence $F(z)$ is
readily given by
\eqn\quadformul{F(z) = {1 \over 2 z^2} \left( 1 -
  \sum_{k = 1}^{\infty}
  {g_k z^{2-k}} - \sqrt{\Delta(z)} \right)}
where $\Delta(z) = \left( 1 - \sum_{k = 1}^{\infty}
  {g_k z^{2-k}} \right)^2 - 4 z^2 \left( 1 - \sum_{k=1}^{\infty} g_k
  \sum_{j=0}^{k-2} z^{2-k+j} F_j \right)$ is the discriminant.
As shown in [\xref\BECA-\xref\BMJ], $\sqrt{\Delta(z)}$ has
a factorization of the form
\eqn\brown{\sqrt{\Delta(z)}= \Gamma(z^{-1}) \sqrt{\kappa(z)}, \qquad
  \Gamma(z^{-1}) = \sum_{q=0}^{\infty} \gamma_q z^{-q}, \qquad
  \kappa(z) = 1 + \kappa_1 z + \kappa_2 z^2}
where the coefficients $\gamma_q$, $\kappa_1$, $\kappa_2$ are power
series in the $g_k$ to be determined. (In the actual proof, one must
for a while assume that degrees are bounded, i.e. that $g_k$ vanishes
for all $k$ larger than some fixed integer, which from \quadformul\
implies that $\sqrt{\Delta(z)}$ does not contain arbitrarily large
negative powers of $z$. In this situation, the existence of a
factorization follows from Brown's theorem \Bro, and $\Gamma(z^{-1})$
is a polynomial in $z^{-1}$.  The bound on degrees can be lifted
eventually and the product $\Gamma(z^{-1}) \sqrt{\kappa(z)}$ has a
well-defined expansion in $z$.)

Let us now explain how to
compute practically $\kappa(z)$ and $\Gamma(z^{-1})$. Following the
approach of \CENSUS, we introduce the power series $R$ and $S$ defined by
\eqn\kappasr{S=-{\kappa_1 \over 2}, \quad
  R={\kappa_1^2 - 4 \kappa_2 \over 16}, \qquad  {\rm i.e.}\quad
  \kappa(z) = (1 - S z)^2 - 4 R z^2.}
This choice is
particularly convenient because the series expansions of
$\sqrt{\kappa(z)}$ and $1/\sqrt{\kappa(z)}$ are related to the
generating functions $P(n;R,S)$ and $P^+(n;R,S)$ for three-step
paths defined in Sect.~1. Namely, it is elementary to check that
\eqn\sqkapexp{\sqrt{\kappa(z)}=1 - S z - 2 R z^2 \sum_{n=0}^{\infty}
  P^+(n;R,S)  z^{n}, \qquad
  {1 \over \sqrt{\kappa(z)}}= \sum_{n=0}^{\infty} P(n;R,S) z^{n}.}
Now, substituting \brown\ in \quadformul, extracting the coefficient 
of $z^n$ for $n \geq 0$
and using the first relation \sqkapexp, we obtain
\eqn\fjtut{F_n= R \sum_{q=0}^{\infty} \gamma_q\, P^+(n+q;R,S)}
which coincides with
\Fpaths\ up to the identification $A_q = \gamma_q R$.  It remains to
determine the unknowns $\gamma_q$, $R$, $S$. For this it is convenient
to rewrite \quadformul\ as
\eqn\gamtut{\sum_{q=0}^{\infty} \gamma_q z^{-q-2} = {
  z^{-2} - \sum\limits_{k=1}^{\infty} g_k z^{-k}
  - 2 F(z) \over \sqrt{\kappa(z)}}}
and recall that $F(z)$ only contains nonnegative powers of $z$, with
$F_0$=1.  Using then the second relation \sqkapexp, extracting the
coefficient of $z^{-q-2}$ for $q \geq 0$ leads directly to \apaths\
while extracting the coefficients of $z^{-1}$ and $z^0$ yields nothing
but \RSeqs. This establishes the general formula for $F_n$ announced
in Sect.~1.

\subsec{Bijective derivation}

\fig{The generating function $F_{n;d}$ for pointed rooted maps with
root degree $n$, where the origin $v$ of the root edge is one of the
vertices closest to the origin $v_0$ of the map among those incident
to the root face, with a distance $d(v_0,v)$ less than $d$, is also that
$Z^+_{d,d}(n)$ for three-step paths of length $n$ starting at, ending
at, and restricted to stay above ordinate $d$, with ordinate-dependent
weights.}{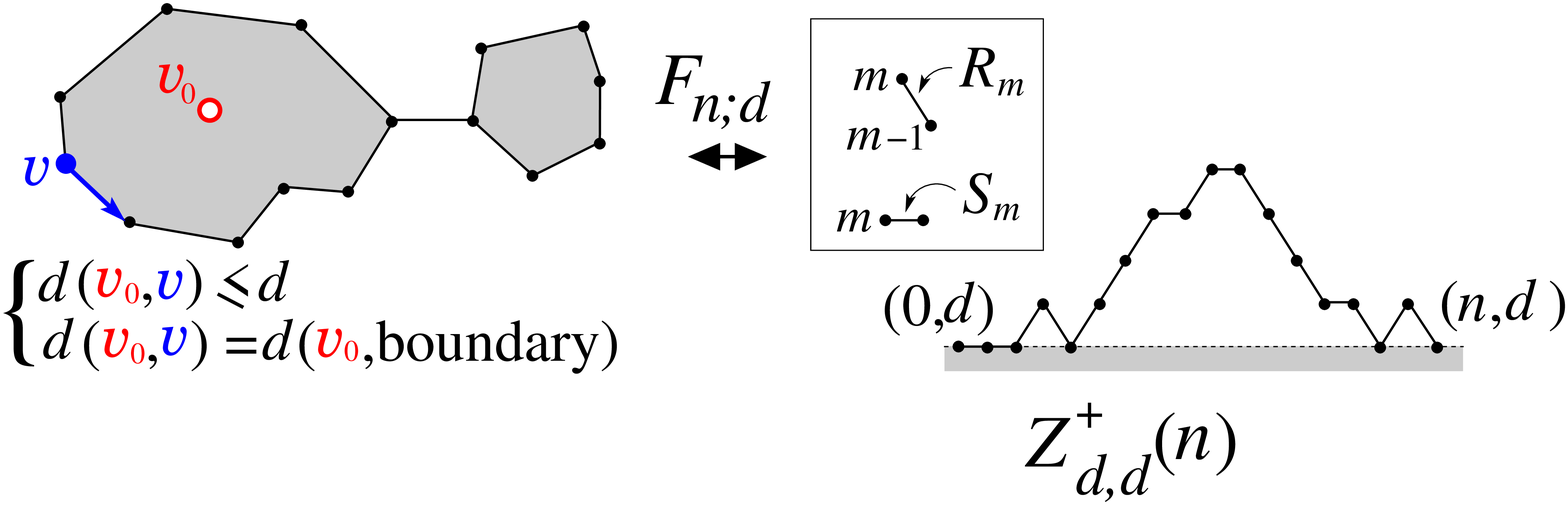}{13.cm} \figlabel\motzmapone

Let us finally present a last derivation of the expression
\Fpaths\ for the moments $F_n$. It is based on a bijective decomposition
of pointed rooted maps keeping track of some distances between the
origin and vertices incident to the root face. Passing to mobiles,
this gives rise to a one-parameter family of expressions for $F_n$
(the so-called {\it conserved quantities}), interesting on their own,
and encompassing \Fpaths-\apaths\ as a limit case.

We first introduce the generating function $f_{n;i}$ ($i\geq 0$) for
pointed rooted maps whose root face has degree $n$ and whose origin is
both at distance $i$ from the origin of the root edge and at distance
$i$ from the {\it boundary} of the root face, defined as the set of
vertices incident to this root face.  In other words, we demand that
the origin of the root edge be one of the vertices closest to the
origin of the map among all vertices incident to the root face. For
$i=0$, the origin of the map is forced to be the origin of the root
edge, so that we have simply $f_{n;0}=F_n$. We also introduce
$F_{n;d}=\sum_{i=0}^d f_{n;i}$ enumerating maps where the above
distance $i$ is less than or equal to some fixed value $d$. Clearly,
we have for any $d\geq 0$ the equality
\eqn\tauto{F_n=F_{n;d}-\sum_{i=1}^d f_{n;i}\ .}
The two terms in the r.h.s may be evaluated in the mobile language as follows:
to evaluate $F_{n;d}$, we look at the sequence of distances from the origin
of the vertices 
around the root face. Starting from the origin of the root edge, we have a 
three-step path from $(0,i)$ to $(n,i)$ for some $i$ between $0$ and $d$, 
with all its ordinates larger than or equal to $i$. This path defines a 
sequence of labels and flags for the labeled and flagged vertices around the 
unlabeled vertex associated with the root face and the full mobile is obtained 
by attaching a rooted mobile to each of these labeled vertices and
a half-mobile to 
each of these flagged vertices. The entire mobile must satisfy (R1) and (R2) 
but, by the same argument as in Sect.~2.3, we may transform it into a mobile 
satisfying (R1) only by 
shifting all labels by $d-i$. The attached rooted mobiles
(respectively half-mobiles) 
are now counted by $R_m$ (respectively $S_m$) with indices corresponding to the
down- (respectively level-) steps of a shifted three-step path from $(0,d)$ 
to $(n,d)$ which never dips below $d$.
To summarize, $F_{n;d}$ is the generating function of three-step paths
from $(0,d)$ to $(n,d)$, staying above $d$, and with a weight $R_m$ 
per down-step $(t,m)\to (t+1,m-1)$ and $S_m$ per level-step $(t,m)\to (t+1,m)$
(see Fig.~\motzmapone\ for an illustration).
With the notations of Sect.~2.1, our first result may be stated as:
\eqn\Fnd{F_{n;d}=Z^+_{d,d}(n)\ .}

\fig{A schematic explanation for the bijection between maps in
$f_{n;i}$ for $i \geq 1$ and maps in ${\cal M}(n;i,k,j)$ with
arbitrary $k \geq 3$ and $j \geq 1$ (see
text).}{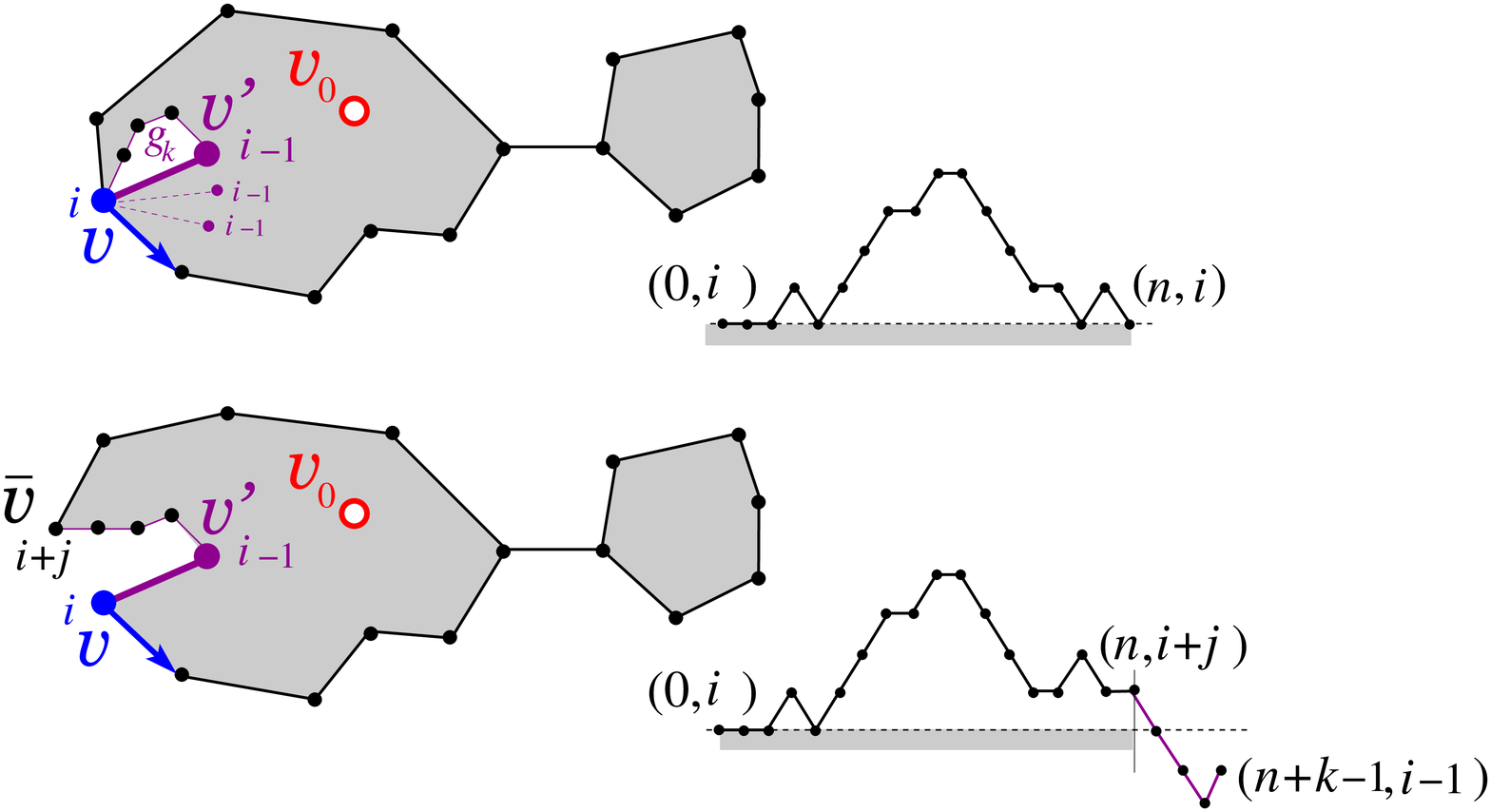}{13.cm} \figlabel\conservproof

As for the subtracted term in \tauto, it concerns maps where the
origin of the root edge, which we from now on denote by $v$, is 
at distance $i\geq 1$ from the origin of the map denoted by $v_0$.
Then there exist
vertices adjacent to $v$ and at distance $i-1$ from $v_0$, hence
edges of the type $(i,i-1)$ incident to $v$. Let us pick the leftmost such 
edge (ordering them from the root edge, see Fig.~\conservproof) and denote 
by $v'$ its endpoint. Note that $v'$ cannot lie on the boundary of the root 
face as, otherwise, the distance from $v_0$ to the boundary would 
be strictly less than $i$. Therefore, the face on the left of the leftmost edge 
from $v$ to $v'$ cannot be the root face but is instead a regular face of 
degree $k\geq 3$ weighted by $g_k$ (the degree cannot be $2$ since 
we picked the leftmost edge from $v$ to $v'$, nor $1$ since it is of
type $(i-1,i)$). We may now merge this face 
with the root face into a new external face with larger degree $n+k$ by 
simply splitting the vertex $v$ into two vertices $v$ and $\bar{v}$ as shown 
in Fig.~\conservproof. Looking again at the sequence of distances along this 
larger boundary starting from the new vertex $v$, the first $n+1$ vertices 
correspond to vertices of the former boundary and their distance to $v_0$ 
is necessarily larger than or equal to what it was before since the splitting 
did not create new adjacences but suppressed some. Moreover, the new vertex 
$v$ is still at distance $i$ from $v_0$ but $\bar{v}$ is necessarily 
at a distance strictly larger than $i$. This is again because we chose the 
leftmost occurrence of an edge $(i,i-1)$ adjacent to $v$. The sequence of 
distances, when going from $v$ to ${\bar v}$ around the new external face, 
now forms a three-step 
path from $(0,i)$ to $(n,i+j)$ for some $j\geq 1$, that never dips below $i$. 
The following distances (going from ${\bar v}$ to $v'$) then form 
a three-step path from $(n,i+j)$ to 
$(n+k-1,i-1)$ with no imposed lower bound (apart from having non-negative 
ordinates), as displayed in Fig.~\conservproof. To summarize, we end up 
with a pointed rooted map, whose root face has degree 
$n+k$ and whose sequence of vertex distances around it defines a three-step 
path made of three parts: a path from $(0,i)$ to $(n,i+j)$ for 
some $j\geq 1$, that never dips below $i$, a three-step path from $(n,i+j)$ 
to $(n+k-1,i-1)$ that never dips below $0$ and a final up-step 
$(n+k-1,i-1)\to (n+k,i)$. We call ${\cal M}(n;i,k,j)$ the family of such maps.

If conversely, we start from a map in ${\cal M}(n;i,k,j)$ for fixed $n \geq 0$
and $i\geq 1$ and arbitrary $k\geq 3$ and $j\geq 1$, gluing the first vertex 
$v$ of the boundary to the $n$-th one $\bar{v}$ into a single vertex 
$v$ will split the external face into a 
face of degree $k$ and a new external face of degree $n$. The new vertex
$v$ will be at distance $i$ from the origin and all the vertices on
the new boundary will be at distance larger than or equal to $i$. Indeed,
the only possible way to reduce the distance from the origin after gluing 
is to use a path which goes through $v$ but such a path has a length larger
than or equal to $i$ so the distance cannot be reduced to less than $i$. 
Finally, all the vertices adjacent to $\bar{v}$ were at a distance
larger than or equal to $i+j-1\geq i$ hence their distance after
gluing remains $\geq i$ so that the leftmost edge $(i,i-1)$ incident 
to $v$ after gluing
is precisely the edge from $v$ to the vertex $v'$ formerly preceding
$v$ along the external face. 
In other words, the splitting and gluing transformations are
inverse of one-another and provide a bijection between maps 
enumerated by $f_{n;i}$ for some fixed $n\geq 0$ and $i\geq 1$ and maps 
in ${\cal M}(n;i,k,j)$ with arbitrary $k\geq 3$ and $j\geq 1$.
These latter maps may themselves be transformed into mobiles
satisfying (R1) and (R2). As before, the set of all values 
$i=1,\ldots, d$ is well captured in a single family of mobiles satisfying 
(R1) only upon shifting the labels and flags by $d-i$. 
In these mobiles, the three-step path ensuring (P) around the unlabeled
vertex associated with the external face is now a path from $(0,d)$ to 
$(n+k,d)$ made of three parts: a first part of length $n$ from $(0,d)$
to $(n,d+j)$ for some $j\geq 1$ that never dips below $d$, enumerated
by $Z^+_{d,d+j}(n)$, a second part of length $k-1$ from $(n,d+j)$ to
$(n+k-1,d-1)$ that never dips below $0$, enumerated by $Z_{d+j,d-1}(k-1)$,
and a final up step from $(n+k-1,d-1)$ to $(n+k,d)$, with weight $1$.  
This is summarized in the equality:
\eqn\exprzero{\sum_{i=1}^d f_{n;i}=\sum_{k\geq 3}g_k\sum_{j\geq 1}
Z^+_{d,d+j}(n)Z_{d+j,d-1}(k-1)\ .}
From \Fnd\ and \exprzero, we deduce eventually the expression
for $F_n$:
\eqn\exprone{F_n=Z^+_{d,d}(n)-\sum_{k\geq 3}g_k\sum_{j\geq 1}
Z^+_{d,d+j}(n)Z_{d+j,d-1}(k-1)}
{\it valid for any} $d\geq 0$.
An important outcome of this result is that the quantity on the r.h.s
of Eq.~\exprone\ is a {\it conserved quantity}, i.e. it is 
{\it independent of} $d$. This conservation property was already proved in 
\DFG\ in the restricted case of bipartite maps via completely different and 
much more involved arguments, with no direct combinatorial interpretation. 
Taking $d=0$, we have in particular the equality $F_n=Z^+_{0,0}(n)$ which 
corresponds precisely to the interpretation of $F_n$ which was discussed 
in Sect.~2.2 and which led us to the continued fraction formula \continued\ 
for $F(z)$. Another interesting limiting case is when $d\to \infty$ in which 
case the r.h.s in \exprone\ may be expressed in terms of $R$ and $S$ only. 
Before we take this limit, let us derive two slightly different, although 
essentially similar, expressions for $F_n$.
They can be obtained by recalling the second recursion relation \rnsnrecur\ for
$R_d$, which be may used to rewrite \exprone\ in the
slightly different form 
\eqn\exprtwo{F_n=Z^+_{d,d}(n) R_d-\sum_{k\geq 2}g_k\sum_{j\geq 0}
Z^+_{d,d+j}(n)Z_{d+j,d-1}(k-1)}
where the sums now include the contribution of $j=0$ and $k=2$ (note that
the term $k=1$ in the first line of \rnsnrecur\ is zero).
The enumerated three-step paths are now paths from $(0,d)$ to $(n+k-1,d-1)$ 
whose first part of length $n$ has ordinates larger than or equal to $d$, with 
the ordinate $d$ now allowed for the abscissa $n$. Denoting by $n+q+1$ the
abscissa of the first occurrence of the ordinate $d-1$, with $q$ between
$0$ and $k-2$, we may write alternatively
\eqn\exprthree{F_n=Z^+_{d,d}(n) R_d-\sum_{k\geq 2}g_k\sum_{q\geq 0}
Z^+_{d,d}(n+q)\, R_d\, Z_{d-1,d-1}(k-2-q)}
valid for all $d\geq 0$.
In particular, by letting $d \to \infty$ and recalling from Sect.~2.1
that $Z_{d,d}(n) \to P(n;R,S)$ and $Z^+_{d,d}(n) \to P^+(n;R,S)$,
we arrive at the simplest formula
\eqn\exprfour{F_n=P^+(n;R,S) R-\sum_{k\geq 2}g_k\sum_{q\geq 0}
P^+(n+q;R,S)\, R\, P(k-2-q;R,S)\ }
which is precisely the desired expression \Fpaths--\apaths.

\fig{The generating function $\sum_{i=0}^{\infty} f_{n;i}$ for pointed
rooted maps with root degree $n$, where the origin $v$ of the root
edge is one of the vertices closest to the origin $v_0$ of the map
among those incident to the root face, is also that $P^+(n;R,S)$ for
Motzkin paths of length $n$ with ordinate-independent
weights.}{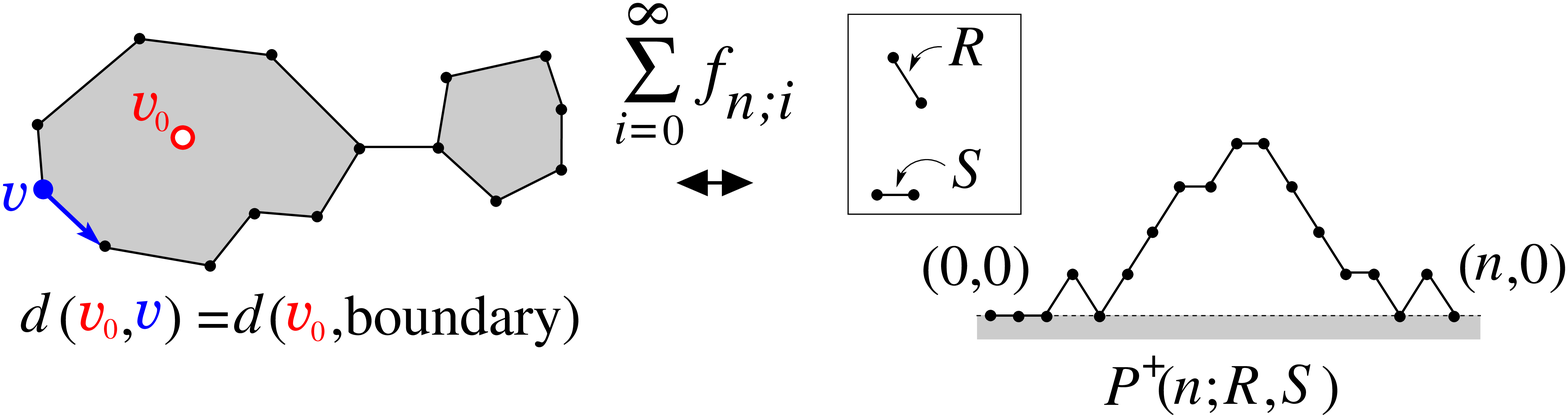}{13.cm} \figlabel\motzmaptwo

A slightly more direct derivation of \exprfour\ consists in using mobiles
without the conditions (R1) and (R2) as done at the very end of Sect.~2.2. 
More precisely, we may in the first place label the vertices of the maps 
counted by 
$f_{n;i}$ by their distance to the origin minus $i$ and consider 
values of $i$ no longer bounded by $d$. 
We end up with mobiles where the conditions
(R1) and (R2) are waived and whose sequence of labels and flags around
the unlabeled vertex associated with the root face match the down- and 
level-steps of a three-step path from $(0,0)$ to $(n,0)$ that never
dips below $0$. In the absence of the conditions (R1) and (R2), each
down- (respectively level-) step is now weighted by $R$ (respectively $S$)
so that we arrive directly at
\eqn\alterfnizero{\sum_{i=0}^\infty f_{n;i}= P^+(n;R,S)\ .}
This formula is illustrated in Fig.~\motzmaptwo.
If we now consider the case of arbitrary $i\geq 1$ and use the above
bijection of Fig.~\conservproof\ with configurations of rooted maps in 
${\cal M}(n;i,k,j)$ for arbitrary $k\geq 3$ and $j\geq 1$, we arrive 
by again shifting the labels in these maps by $i$ to mobiles counted 
directly by
\eqn\alterfnione{\sum_{i=1}^\infty f_{n;i}=
\sum_{k\geq 3} g_k \sum_{j\geq 1} \sqrt{R}\, P^+_j(n;R,S) P_{-j-1}(k-1;R,S)}
where, for an arbitrary integer $m$, $P_m(q;R,S)$ denotes as before the 
generating function for arbitrary three-step paths from $(0,0)$ to $(q,m)$
with a weight $\sqrt{R}$ per up- or down-step and $S$ per level-step and,
for $m\geq 0$, $P^+_m(n;R,S)$ restricts the enumeration to paths whose 
ordinates all stay above $0$.
Using the identity \rsrecurbis\ for $R$ we arrive directly at 
\eqn\exprFnbis{\eqalign{F_n& =P^+(n;R,S) R-\sum_{k\geq 2}g_k\sum_{j\geq 0}
\sqrt{R} P^+_{j}(n;R,S)P_{-j-1}(k-1;R,S) \cr &=
P^+(n;R,S) R-\sum_{k\geq 2}g_k\sum_{q\geq 0}
P^+(n+q;R,S)\, R\, P(k-2-q;R,S) \cr}}
upon denoting, in the r.h.s of  the second line, by $n+q+1$ the abscissa of 
the first occurrence of the ordinate $-1$ in the paths enumerated 
in the r.h.s of the first line. 

\newsec{Study of the Hankel determinants}

Having established the formula \Fpaths\ for $F_n$, we now study its
consequences on the Hankel determinants $\displaystyle{H_n = \det_{0
\leq i,j \leq n} F_{i+j}}$ and minors $\displaystyle{{\tilde H}_n =
\det_{0 \leq i,j \leq n} F_{i+j+\delta_{j,n}}}$.  We first establish
formulas \Hdetb\ and \tHdetb\ via a natural path decomposition
which translates into matrix identities (Sect.~4.1). We then identify the
resulting determinants with symplectic Schur functions (Sect.~4.2).

\subsec{Path decomposition}

\fig{A schematic picture of the path decomposition leading to
(4.1).}{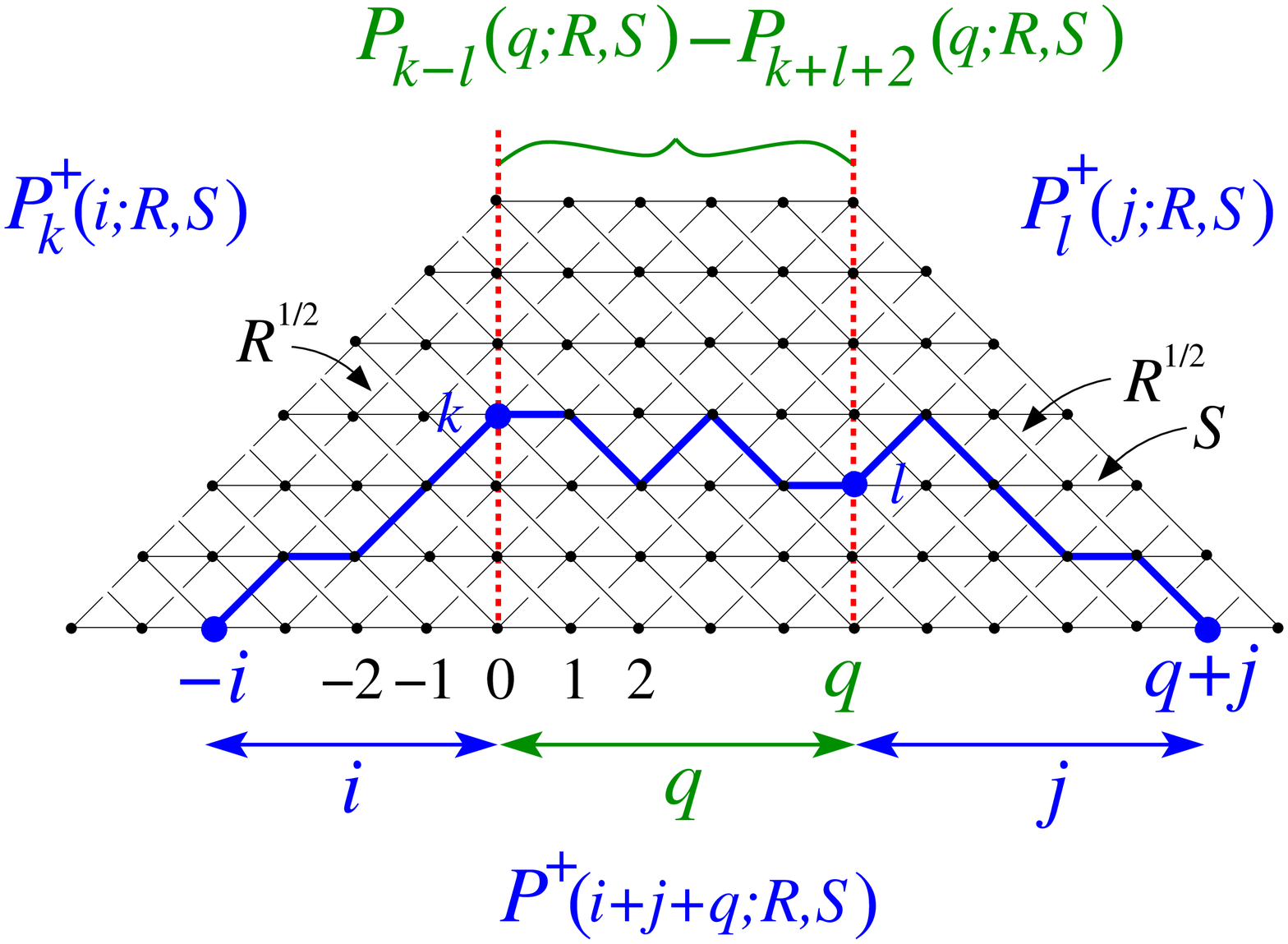}{11.cm} \figlabel\decomp

We first establish the formulas \Hdetb\ and \tHdetb. We start with the identity
\eqn\decompkl{P^+(i+j+q;R,S)=\sum_{k=0}^i \sum_{\ell=0}^j P^+_k(i;R,S)\,
\left(P_{k-\ell}\left(q;R,S\right)-P_{k+\ell+2}\left(q;R,S\right)\right)\,
P^+_{\ell}(j;R,S)}
obtained by decomposing the paths of length $i+j+q$ in $P^+(i+j+q;R,S)$ into
a first path of length $i$ ending at the ordinate $k$, a second path of 
length $q$
ending at the ordinate $\ell$ and a third path of length $j$ back to 
ordinate $0$ (see
Fig.~\decomp\ for an illustration). The medial part must have non-negative
ordinates only, hence we must subtract from $P_{k-\ell}(q;R,S)$ the
configurations which dip below $0$. By a classical reflection argument, those
configurations are in one-to-one correspondence with paths from ordinate
$-k-2$ to ordinate $\ell$, enumerated by $P_{k+\ell+2}(q;R,S)$.

Using now formula \Fpaths, we deduce the matrix identity:
\eqn\matrixid{{\bf H}={\bf T}^{\rm T} \cdot {\bf B} \cdot {\bf T}}
where ${\bf H}$, ${\bf T}$ and ${\bf B}$ are the semi-infinite matrices
\eqn\matrices{
{\bf H} \equiv (F_{i+j})_{i,j\geq 0,}\quad
{\bf T} \equiv (P^+_{i}(j;R,S))_{i,j\geq 0},\quad
{\bf B} \equiv (B_{i-j}-B_{i+j+2})_{i,j\geq 0},}
with $B_i$ defined as in \Hdetb\ and ${\bf T}^{\rm T}$ denoting the transpose of
${\bf T}$. Note that ${\bf T}$ is upper triangular, hence the matrix product
\matrixid\ involves only a finite sum for each matrix element of ${\bf H}$.

For a semi-infinite matrix ${\bf M}=(M_{i,j})_{i,j \geq 0}$, let us
denote by ${\bf M}_n$ its submatrix obtained by keeping only the row
and column indices between $0$ and $n$. Note that ${\bf M}_n$ has
size $(n+1) \times (n+1)$. We also denote by ${\bf M}_n^{(i;)}$,
${\bf M}_n^{(;j)}$ and ${\bf M}_n^{(i;j)}$ the submatrices obtained by
further removing respectively the row with index $i$, the column with
index $j$ and both row $i$ and column $j$.

We have, by definition, $H_n = \det {\bf H}_n$. 
Because ${\bf T}$ is upper triangular,
the identity \matrixid\ restricts safely to
\eqn\matrixidn{{\bf H}_n={\bf T}_n^{\rm T} \cdot {\bf B}_n \cdot {\bf T}_n}
which yields
\eqn\detid{H_n = \left( \det {\bf T}_n \right)^2
  \det {\bf B}_n.}
Noting that $T_{ii}=P^+_i(i)=R^{i/2}$, the triangular determinant is
readily computed, establishing formula \Hdetb.

As for $\tilde{H}_n=\det {\bf H}_{n+1}^{(n+1;n)}$, we have the restriction
\eqn\matrixidmin{{\bf H}_{n+1}^{(n+1;n)}={\bf T}_n^{\rm T} \cdot
  {\bf B}_{n+1}^{(n+1;)} \cdot {\bf T}_{n+1}^{(;n)}}
and the Cauchy-Binet formula yields
\eqn\detcb{\tilde{H}_n= \det {\bf T}_n \sum_{i=0}^{n+1}
  \det {\bf B}_{n+1}^{(n+1;i)} \det {\bf T}_{n+1}^{(i;n)}.}
Note that ${\bf T}_{n+1}^{(i;n)}$ is still triangular and, unless
$i=n$ or $n+1$, it has at least one zero on the diagonal. Therefore only those
two terms are non-zero in \detcb\ and, as $T_{n,n+1}=P^+_n(n+1)=(n+1) R^{n/2} S$, the triangular determinants are again readily computed. Finally we have
$\det {\bf B}_{n+1}^{(n+1;n+1)}=\det {\bf B}_n=H_n/R^{n(n+1)/2}$, and
$\det {\bf B}_{n+1}^{(n+1;n)}$ is nothing but the determinant in formula
\tHdetb, which is now established.

To conclude this subsection, let us mention that using the
decomposition \matrixid\ of the Hankel matrix ${\bf H}$ for
determinant evaluation amounts to doing elementary row and column
manipulations, as stated in Sect.~1.

\subsec{Determinants and minors of the matrix ${\bf B}$ identified as Schur functions}

We now turn to the study of the quantities $\det {\bf B}_n$ and $\det
{\bf B}_n^{(n+1;n)}$ encountered above, which are minors of the
semi-infinite matrix ${\bf B}$. This matrix is remarkably the
difference between the symmetric Toeplitz matrix $(B_{i-j})_{i,j \geq
0}$ and the Hankel matrix $(B_{i+j+2})_{i,j \geq 0}$, strongly
suggesting possible determinantal identities. The Laurent series
$\sum_{i=-\infty}^{\infty} B_i x^i$ is known as the symbol of the
Toeplitz matrix.

To proceed further with our analysis, we shall assume that the map
degrees are bounded, namely we fix a positive integer $p$ and consider
maps whose faces have degree at most $p+2$. As already mentioned
above, this amounts to setting $g_k=0$ for $k>p+2$, so that we only
deal with power series in the finite set of variables
$g_1,\ldots,g_{p+2}$. Then $B_i$ vanishes for $\vert i \vert >p$ (and
$B_p \neq 0$) so that ${\bf B}$ is a band matrix. The symbol is a
symmetric Laurent polynomial of degree $p$ in $x$, which thus admits
$2p$ roots of the form ${\bf
x}=(x_1,1/x_1,x_2,1/x_2,\ldots,x_p,1/x_p)$. These are precisely the
solutions of Eq.~\chareqb\ (which by analogy with the
``perturbative'' approach of \GEOD\ we call the characteristic
equation). Up to the overall factor $B_p=-g_{p+2}R^{(p+2)/2}$, 
we may view ${\bf x}$ as a
parametrization of the $B$'s, namely
\eqn\Btoe{B_j = (-1)^{p+j} B_p\, e_{p+j}({\bf x})}
where $e_i({\bf x})$ is the $i$th elementary symmetric function of ${\bf x}$,
defined for instance via
\eqn\defe{\prod_{i=1}^p(1+x_i t)\left(1+{t\over x_i}\right)=
\sum_{i=-\infty}^\infty
e_i({\bf x}) t^i\ .}
Note that $e_{p-j}({\bf x})=e_{p+j}({\bf x})$ and $e_i({\bf x})=0$ for 
$i>2p$ or $i<0$.
It follows that $\det {\bf B}_n$ and $\det
{\bf B}_n^{(n+1;n)}$ are symmetric functions of ${\bf x}$, which we may
hopefully identify.

Let us first quickly describe a naive strategy, based on finding the
zeroes of the symmetric functions at hand. We observe that, for $x$ a
root of the symbol, the semi-infinite vector $\left(
(x^{i+1}-x^{-i-1})/(x-x^{-1}) \right)_{i \geq 0}$ is in the kernel of
${\bf B}$. Its restriction to the indices $\{0,\ldots,n\}$ is however
generally not in the kernel of ${\bf B}_n$ nor ${\bf B}_n^{(n+1;n)}$.
We may attempt to find a vector in the kernel of ${\bf B}_n$ (resp.\
${\bf B}_n^{(n+1;n)}$) by taking a linear combination over the $p$
independent roots of the symbol. Doing so, we find that we have to satisfy
$p$ ``boundary conditions'' which are linear equations for the $p$
coefficients of the linear combination. A non-trivial solution exists
if the associated $p \times p$ determinant vanishes, hence this
determinant (which is again a symmetric function of ${\bf x}$) divides
$\det {\bf B}_n$ (resp.\ $\det {\bf B}_n^{(n+1;n)}$). By comparing
degrees, we find that they only differ by a constant. In the end, we
have expressions for $\det {\bf B}_n$ and $\det {\bf B}_n^{(n+1;n)}$
in terms of $p \times p$ determinants. It turns out that these are
nothing but the determinants present in the Weyl character formula for
the symplectic group ${\rm Sp}_{2p}$, written down below.
What we have just encountered are classical formulas in
representation theory.

Let us now be more educated and directly identify $\det {\bf B}_n$ and
$\det {\bf B}_n^{(n+1;n)}$ as instances of the general ``symplectic
$e$-formula'' [\xref\FULHA,\xref\FULKRA]
\eqn\eformula{{\rm sp}_{2p}(\lambda,{\bf x}) = \det_{1 \leq i,j \leq m}
  \left( e_{\lambda'_j -j+i}({\bf x}) - e_{\lambda'_j-j-i}({\bf x}) \right).}
Here ${\rm sp}_{2p}(\lambda,{\bf x})$ stands for the {\it symplectic
Schur function} associated with the partition $\lambda$ (having at most $p$
parts), $\lambda'$ denoting its conjugate partition (having at most $m$
parts). Comparing with \Btoe\ and recalling that $e_{p+q}({\bf
x})=e_{p-q}({\bf x})$ for all $q$, we obtain the general formula
\eqn\minorsB{ \det_{0 \leq i,j \leq n} (B_{i-j-\mu_j} - B_{i+j+\mu_j+2}) =
  (-1)^{|\lambda|} B_p^{n+1} {\rm sp}_{2p}(\lambda,{\bf x})}
where $\mu_j=p-\lambda'_{j+1}$, $m=n+1$ and $|\lambda|$ denotes the
sum of $\lambda$ (equal to that of $\lambda'$). It is now immediate to
identify $\det {\bf B}_n$ with the case $\mu_j=0$, hence with the
``rectangular'' partition $\lambda'_{p,n+1}=p^{n+1}$ made of $n+1$ parts of
size $p$, conjugate to the partition $\lambda_{p,n+1}=(n+1)^p$ with $p$ parts
of size $n+1$. Similarly $\det {\bf B}_n^{(n+1;n)}$ corresponds to the
case $\mu_j=\delta_{j,n}$, hence to the ``nearly-rectangular''
partition $\tilde{\lambda}'_{p,n+1}=p^n(p-1)$ made of $n$ parts of size 
$p$ and one part of size $p-1$, conjugate to the partition
$\tilde{\lambda}_{p,n+1}=(n+1)^{p-1}n$ with $p-1$ parts of size $n+1$ 
and one part of size $n$. 
We end up with the compact formula for the Hankel determinants
\eqn\Hnsp{H_n= 
(-1)^{p(n+1)} B_p^{n+1}\, 
R^{{n(n+1)\over 2}}
\ {\rm sp}_{2p}(\lambda_{p,n+1},{\bf x})} 
from which, using \Rndet, we deduce
\eqn\RSchur{R_n = R \, {{\rm sp}_{2p}(\lambda_{p,n+1},{\bf x}) \,
  {\rm sp}_{2p}(\lambda_{p,n-1},{\bf x}) \over
  {\rm sp}_{2p}(\lambda_{p,n},{\bf x})^2}\ .}
Similarly, for the Hankel minors we have
\eqn\tHnsp{\tilde{H}_n-(n+1)S\, H_n= 
(-1)^{p(n+1)+1} B_p^{n+1}\, R^{{n^2+n+1\over 2}} 
\ {\rm sp}_{2p}(\tilde{\lambda}_{p,n+1},{\bf x})}
from which, using \Sndet, we deduce
\eqn\SSchur{S_n = S - \sqrt{R} \left(
 {{\rm sp}_{2p}(\tilde{\lambda}_{p,n+1},{\bf x}) \over 
  {\rm sp}_{2p}(\lambda_{p,n+1},{\bf x}) } -
 {{\rm sp}_{2p}(\tilde{\lambda}_{p,n},{\bf x}) \over 
  {\rm sp}_{2p}(\lambda_{p,n},{\bf x}) } \right)\ .}

Beside the above $e$-formula, other expressions are known for the
symplectic Schur function ${\rm sp}_{2p}$, namely the ``$h$-formula''
and the Weyl character formula. Both involve determinants of size $p$,
independently of the variable $n$ in
$H_n$, $\tilde{H}_n$, $R_n$ and $S_n$.
On the one hand, the $h$-formula involves the complete symmetric 
function $h_i({\bf x})$, 
defined for instance via
\eqn\defh{\prod_{i=1}^p{1\over 1-x_i t}{1\over 1-{t\over x_i}}= 
\sum_{i=0}^\infty h_i({\bf x}) t^i\ ,}
and reads [\xref\FULHA-\xref\FULKRA]
\eqn\hformula{{\rm sp}_{2p}(\lambda,{\bf x}) = 
\det_{1 \leq i,j \leq p}
  \left(h_{\lambda_j-j+1}({\bf x}) \, \raise -2pt \vdots\,  
h_{\lambda_j -j+i}({\bf x}) + 
h_{\lambda_j-j-i+2}({\bf x}) \right)}
where $(a_j\, \raise -2pt \vdots\, a_{i,j})$ denotes the matrix with 
elements $a_j$ in the first row and elements $a_{i,j}$ in the rows $i>1$. 
It immediately allows to rewrite the r.h.s of \Hnsp\ and \tHnsp\ as
determinants of size $p$. On the other hand, 
the Weyl character formula reads \FULHA
\eqn\Weylformu{{\rm sp}_{2p}(\lambda,{\bf x}) =
  { \det\limits_{1 \leq i,j \leq p} (x_i^{\lambda_{p+1-j}+j} - 
x_i^{-\lambda_{p+1-j}-j}) \over
    \det\limits_{1 \leq i,j \leq p} (x_i^{j} - x_i^{-j}) }}
and involves the ratio of two determinants of size $p$. However, the
denominators cancel in the expressions \RSchur\ for $R_n$ and \SSchur\ 
for $S_n$, which establishes the announced nice 
formulas \finalRn\ and \Snfinal.

\newsec{The special case of maps with even face degrees}

As mentioned in Sect.~1, some simplifications occur when we restrict
our analysis to bipartite maps. The distance from a given origin
in the map changes parity between adjacent vertices and consequently, 
there are no edges of type $(n,n)$, leading to mobiles {\it without flagged 
vertices}. At the level of generating function, this implies that $S$ and
$S_n$ for all $n$ vanish so that the continued fraction \continued\ 
becomes of the Stieltjes type \Stieltjes, consistent with the property that
$F_n=0$ for $n$ odd as the root face must have even degree in a bipartite map.
The Hankel determinants factorize as 
\eqn\factorhankbis{H_{2n}= h^{(0)}_n h^{(1)}_{n-1}, \qquad
H_{2n+1}= h^{(0)}_n h^{(1)}_{n}}
where 
\eqn\hzeroone{h^{(0)}_n=\det_{0\leq i,j\leq n} F_{2i+2j}, \qquad
h^{(1)}_n=\det_{0\leq i,j\leq n} F_{2i+2j+2}\ .
}
Using the same path decomposition as in Sect.~4.1, we may write 
\eqn\hzerooneb{\eqalign{
h^{(0)}_n& =  R^{n(n+1)} \det_{0\leq i, j \leq n} ({\hat B}_{i-j}
-{\hat B}_{i+j+1}) \cr
h^{(1)}_n& =  R^{(n+1)^2} \det_{0\leq i, j \leq n} ({\hat B}_{i-j}
-{\hat B}_{i+j+2})\ . \cr}}
where we introduce the notation
\eqn\hatb{{\hat B}_i= B_{2i}}
for all integers $i$ (note that $B_{2i+1}=0$ and ${\hat B}_{-i}={\hat B}_i$).
Let us now turn to the case of maps with maximal degree $2p+2$, for 
some $p\geq 1$, so that ${\hat B}_i$ vanishes for $\vert i \vert >p$. 
Writing the characteristic equation as
\eqn\charpair{\sum_{i=-p}^{p} {\hat B}_i\, y^i =0, \qquad y\equiv x^2, }
whose $2p$ roots are gathered in the $2p$-uple 
${\bf y}\equiv (y_1,1/y_1, \ldots, y_p,1/y_p)$,
we now have the identification
\eqn\hatbe{{\hat B}_i= (-1)^{p+i}\, {\hat B}_p\, e_{p+j}({\bf y})}
in terms of the elementary symmetric function $e_{j}({\bf y})$.
This yields 
\eqn\hzerooneSchur{\eqalign{
h^{(0)}_n & =(-1)^{p(n+1)}\, ({\hat B}_p)^{n+1}\, R^{n(n+1)}\ 
\det_{1\leq i,j \leq n+1}\left(e_{p-j+i}({\bf y})+ e_{p-i-j+1}({\bf y})\right) 
\cr 
&= (-1)^{p(n+1)}\, ({\hat B}_p)^{n+1}\, R^{n(n+1)}\
{\rm o}_{2p+1}(\lambda_{p,n+1},{\bf y})\cr
h^{(1)}_n & =(-1)^{p(n+1)}\, ({\hat B}_p)^{n+1}\, R^{(n+1)^2}\ 
\det_{1\leq i,j \leq n+1}\left(e_{p-j+i}({\bf y})- e_{p-i-j}({\bf y})\right) 
\cr  &= 
(-1)^{p(n+1)}\, ({\hat B}_p)^{n+1}\, R^{(n+1)^2}\
{\rm sp}_{2p}(\lambda_{p,n+1},{\bf y})\cr
}}
where, as before, $\lambda_{p,n+1}=(n+1)^p$ is the partition with $p$ parts
of size $n+1$.  Here we used again \eformula\ to identify the second 
determinant with a symplectic Schur function of ${\bf y}$, while the
first determinant is now recognized as an instance of the general
formula
\eqn\eformulaortho{\eqalign{{\rm o}_{2p+1}(\lambda,{\bf y}) & = 
\det_{1 \leq i,j \leq m} \left( e_{\lambda'_j -j+1}({\bf {\tilde y}})
\, \raise-2pt \vdots \, 
e_{\lambda'_j -j+i}({\bf {\tilde y}}) + 
e_{\lambda'_j-j-i+2}({\bf {\tilde y}}) \right) 
\cr & = \det_{1 \leq i,j \leq m} \left( e_{\lambda'_j -j+i}({\bf y}) + 
e_{\lambda'_j-j-i+1}({\bf y}) \right) 
\cr}
}
where ${\rm o}_{2p+1}(\lambda,{\bf y})$ stands for the {\it odd-orthogonal
Schur function} associated with the partition $\lambda$ (having at most $p$ 
parts), $\lambda'$ denoting its conjugate partition (having at most $m$
parts), and
where ${\bf {\tilde y}}$ denotes the $(2p+1)$-uple
$(y_1,1/y_1, \ldots, y_p,1/y_p,1)$ with an additional 
$1$ term. The first identity may be found in \FULKRA\ and it implies 
the second one by elementary manipulations on the determinant.
Again, the odd-orthogonal Schur function admits a simple Weyl formula, 
namely:
\eqn\Weylformuortho{{\rm o}_{2p+1}(\lambda,{\bf y}) =
  { \det\limits_{1 \leq i,j \leq p} (y_i^{\lambda_{p+1-j}+j-{1\over 2}} - 
y_i^{-\lambda_{p+1-j}-j+{1\over 2}}) \over
    \det\limits_{1 \leq i,j \leq p} (y_i^{j-{1\over 2}} - y_i^{-j+
{1\over 2}}) }\ .}
Alternatively, the factorization \factorhankbis\ of $H_n$, with $h^{(0)}_n$ 
and $h^{(1)}_n$ expressed directly as in \hzerooneSchur, may be obtained
right away from the general formula \Hnsp\ (with $p$ replaced by $2p$)
for $H_n$ upon using the identity
\eqn\Schurident{{\rm sp}_{4p}(\lambda_{2p,n+1}; {\bf x}\cup {\bf -x})
= (-1)^{p(n+1)}\, {\rm sp}_{2p}(\lambda_{p, \left\lfloor {n+1\over 2}
\right\rfloor},{\bf y})\, {\rm o}_{2p+1}(\lambda_{p, \left\lfloor {n+2\over 2}
\right\rfloor },{\bf y})}
where ${\bf x}\cup {\bf -x}=(x_1,1/x_1,\ldots,x_p,1/x_p,-x_1,-1/x_1,\ldots
-x_p,-1/x_p)$ and where 
${\bf y}=(y_1,1/y_1,\ldots, y_p,1/y_p)$ with $y_i=(x_i)^2$.
This identity is easily proved via elementary determinant manipulations
in the Weyl formula for each of these Schur functions.

The final ``nice'' formula for $R_n$ becomes 
\eqn\finalRnpair{R_n=R \ 
{\det\limits_{1\leq i,j \leq p} \left(y_i^{{n\over 2}+j+{1\over 2}}
-y_i^{-{n\over 2}-j-{1\over 2}}\right) 
\, \det\limits_{1\leq i,j \leq p} \left(y_i^{{n\over 2}+j-1}
-y_i^{-{n\over 2}-j+1}\right)  \over 
\det\limits_{1\leq i,j \leq p} \left(y_i^{{n\over 2}+j}
-y_i^{-{n\over 2}-j}\right) \, 
\det\limits_{1\leq i,j \leq p} \left(y_i^{{n\over 2}+j-{1\over 2}}
-y_i^{-{n\over 2}-j+{1\over 2}}\right)} 
}
which matches formula (5.8) of Ref.~\GEOD, given there without a
complete proof. 
To make the identification
complete, we have to rewrite the characteristic equation \chareqP\ as 
\eqn\charpairbis{\eqalign{1& =
\sum_{k=1}^{p+1} g_{2k} \sum_{q=0}^{k-1} P(2k-2-2q;R,0)
\left(\sqrt{R}\, x+{\sqrt{R}\over x}\right)^{2q}\cr
&= 
\sum_{k=1}^{p+1} g_{2k} \, R^{k-1}\, \sum_{q=0}^{k-1} {2k-2-2q \choose
k-q-1} \left(x+{1\over x}\right)^{2q} \cr &=
\sum_{k=1}^{p+1} g_{2k} \, R^{k-1}\, \sum_{m=0}^{k-1} {2k-1 \choose
k-m-1} \sum_{j=-m}^m y^j 
\ ,\cr}}
with $y=x^2$. This is precisely the form found in \GEOD.
The third step in \charpairbis\ follows from the identity
\eqn\idbinom{
  \sum_{q=|j|}^{k-1} {2k-2q-2 \choose k-q-1} {2 q\choose q+j} =
  \sum_{m=|j|}^{k-1} {2k-1\choose k-m-1}}
which itself follows from the more obvious identity
\eqn\idbinombis{
  \sum_{q=|j|}^{k-1} {2k-2q-2\choose k-q-1} \left[{2q\choose q+|j|} - 
{2 q\choose q+|j|+1} \right] = {2k-1\choose k-|j|-1} }
obtained for instance by enumerating paths with $\pm 1$ steps of length
$2k-1$ starting at height $0$ and ending at height $2|j|+1$ in two
manners: either directly (r.h.s) or by decomposing at the last passage
at height $0$ (l.h.s).

\newsec{Triangulations and quadrangulations}

We now turn to the particularly simple cases of triangulations and
quadrangulations, where the expressions \finalRn, \Snfinal\ (for
triangulations) and \finalRnpair\ (for quadrangulations) involve only
determinants of size $1$. They admit an elementary combinatorial
interpretation in terms of one-dimensional hard dimers, as we shall
see.

\subsec{Triangulations}

\fig{Interpretation of formula (6.3) for $F_{i+j}$ in the case of
  triangulations. The weights per step are different in the central
  strip as shown.}{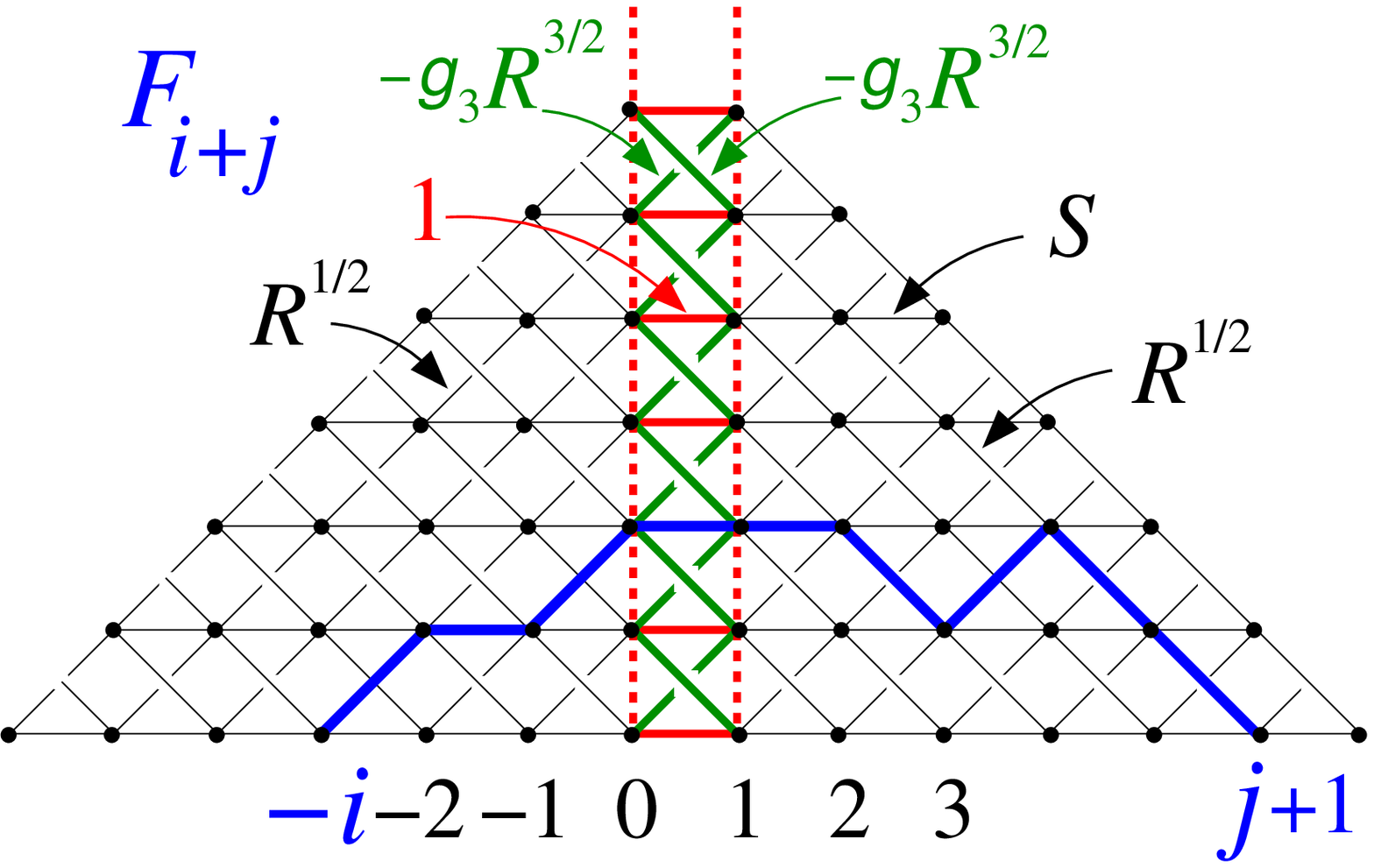}{8.cm} \figlabel\triang
\fig{Correspondence between configurations of non-intersecting paths
on the graph of Fig.~\triang\ and hard dimers on a segment of length
$n$.}{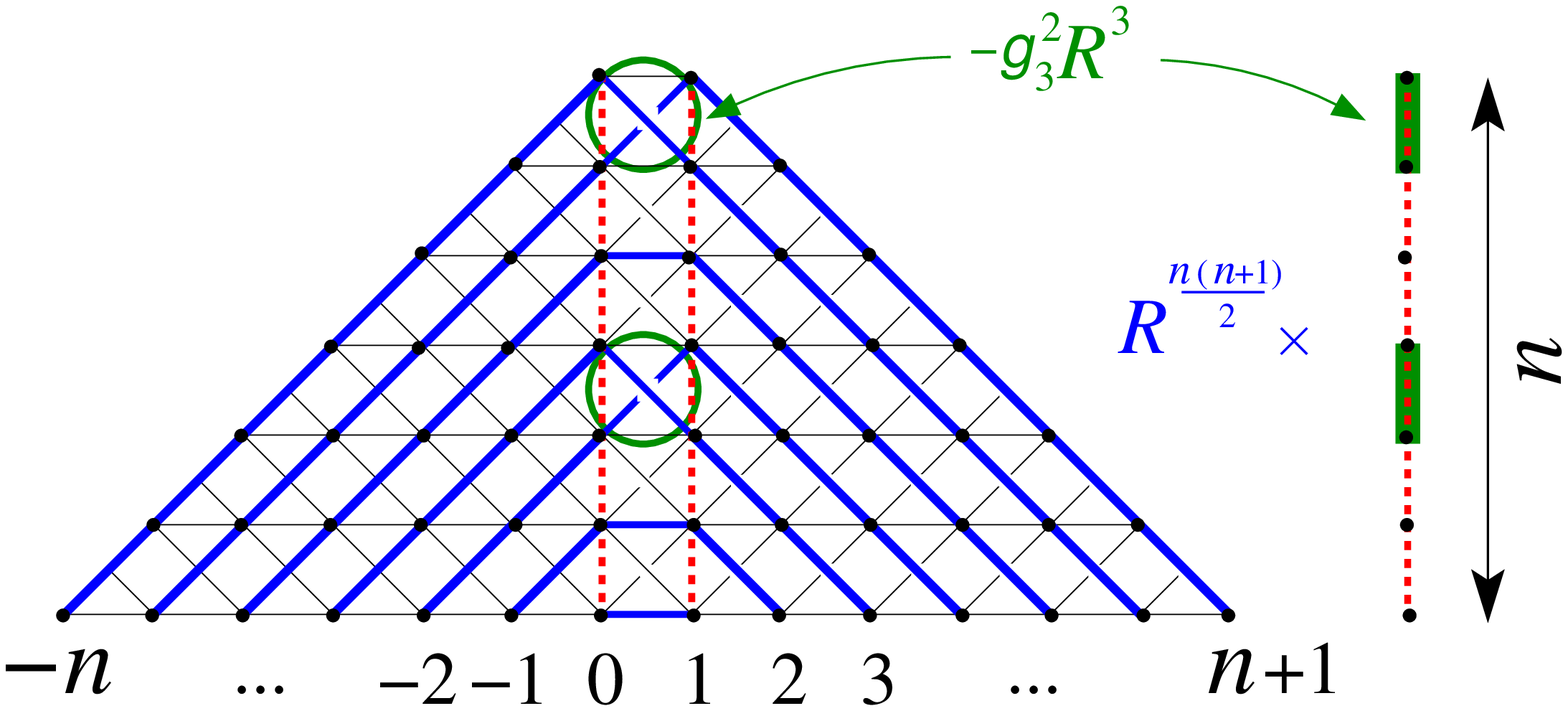}{10.cm} \figlabel\triangHO
The case of triangulations, characterized by $g_k=g_3 \delta_{k,3}$ 
is particularly simple. In this case, the only non-vanishing $A_q$ coefficients
are $A_0$ and $A_1$, given from \apaths\ by
\eqn\azeroone{A_0=R-g_3\, R\, S, \qquad A_1=-g_3\, R}
with $R$ and $S$ given implicitly from \RSeqs\ by 
\eqn\RStrian{S=g_3(S^2+2 R), \qquad R=1+{g_3\over 2}(S^3+6 g_3\, R\, S)
-{S^2\over 2}= 1+ 2 g_3\, R\, S\ .}
From \Hdetb, the non-vanishing $B_i$ coefficients are
\eqn\Btriang{\eqalign{B_0&= A_0+A_1\, S =1-2 g_3\, R\, S=1\cr
B_1&= B_{-1}= A_1\, \sqrt{R}=- g_3\, R^{3/2}\cr}}
while the characteristic equation \chareqb\ reads
\eqn\chartriang{1-g_3\, R^{3/2}\left(x+{1\over x}\right)=0\ .}
Formulas \finalRn\ and \Snfinal\ for the two-point functions reduce to
\eqn\finaltri{\eqalign{
R_n &= R {(x^{n+2}-x^{-n-2})(x^{n}-x^{-n})\over (x^{n+1}-x^{-n-1})^2} \cr
S_n &= S- \sqrt{R}\left({x^{n+1}-x^{-n-1} \over x^{n+2}-x^{-n-2}}-
{x^{n}-x^{-n} \over x^{n+1}-x^{-n-1}} \right)\ .\cr}}
To recover the expressions given in \PDFRaman, we set
\eqn\xy{y=x^2\ ,} 
which is solution of 
\eqn\ytriang{y+{1\over y}+2={1\over g_3^2\, R^3}\ ,}
and obtain eventually
\eqn\RnSntrian{\eqalign{R_n&=R {(1-y^n)(1-y^{n+2})\over (1-y^{n+1})^2}\cr
S_n& =S-g_3 R^2 y^n {(1-y)(1-y^2)\over (1-y^{n+1})(1-y^{n+2})}\ .\cr}}

Interestingly, these formulas admit a combinatorial interpretation 
in terms of one-dimensional hard-dimers. Indeed, for triangulations,
formula \Fpaths\ for the generating function $F_n$ reduces to 
\eqn\Fntrian{F_n= A_0\, P^+(n;R,S)+A_1\, P^+(n+1;R,S)}
so that $F_{i+j}$ may be interpreted as enumerating three-step paths with 
non-negative ordinates from, say, $(-i,0)$ to $(j+1,0)$ with a weight 
$\sqrt{R}$ per up- or down-step and $S$ per level-step, except for steps
in a central strip of width $1$ (i.e. the strip between abscissas $0$ and $1$)
which receive instead weights incorporating the $A_q$ factors, namely 
$A_1 \times \sqrt{R}=-g_3\, R^{3/2}$ for the up- and down-steps,
and $A_0\times 1 +A_1\times S= R-2 g_3\, R\, S=1$ for the level-steps 
(see Fig.~\triang). Note that these weights are nothing but $B_1$ and $B_0$ 
respectively, as might be seen by specializing the discussion of Sect.~4.1.
Alternatively, the paths may be viewed as oriented
paths on a graph drawn from a square grid in the upper-half plane by keeping 
the horizontal sides and the diagonals of the squares, with the two 
diagonals in each square viewed as 
passing on top of each other with no vertex at their crossing point
in the plane. The graph is implicitly endowed with some orientation from left
to right so that the paths always have strictly increasing abscissas. On such 
an oriented graph, we may apply the so-called LGV (Lindstr\"om-Gessel-Viennot) 
lemma \LGV, which states that 
$H_n=\det_{0\leq i,j\leq n}F_{i+j}$ enumerates configurations made of $n+1$ 
paths connecting the set of points $(I_j)_{j=0,\cdots n}$ with coordinates 
$(-j,0)$ to the set of points $(O_j)_{j=0,\cdots,n}$ with coordinates 
$(j+1,0)$ {\it with no intersections at vertices of the graph}. In addition
to the above weights for up-, down- and level-steps, each configuration 
receives a $\pm 1$ factor equal to the signature of the permutation $\sigma$ 
of $\{0, \cdots, n\}$ characterizing its connections (namely $I_j$ is
connected to $O_{\sigma(j)}$). Now it is easy to see that, for such
configurations of non-intersecting paths, the part of the paths outside the 
central strip is entirely fixed to be
a set of straight lines going up on the left of the strip and straight lines
going down on its right, contributing an overall factor $R^{n(n+1)/2}$ 
(see Fig.~\triangHO). The only freedom comes from the possibility of crossings 
(in the plane) along the two diagonals of a given square in 
the strip (since these diagonals do not intersect on the graph). Moreover,
two such crossings cannot take place on two neighboring squares in the
strip. An acceptable crossing configuration may therefore be viewed 
as a {\it hard dimer} configuration on the segment $[0,n]$, as shown in
Fig.~\triangHO. Each dimer receives a weight 
$W=-(-g_3\, R^{3/2})^2=-g_3^2\, R^3$ with a minus sign for each crossing,
equal to the signature of a transposition component in the permutation
$\sigma$. 
In conclusion, the Hankel determinant $H_n$ is simply given by
\eqn\Hntrian{H_n=R^{{n(n+1)\over 2}}Z_{{\rm hard\ dimers}\atop {\rm on\ [0,}n
{\rm ]}}}
where $Z_{{\rm hard\ dimers}\atop {\rm on\ [0,}n
{\rm ]}}$ is the generating function for hard dimers on 
the segment $[0,n]$ with weight $W$ per dimer. Now it is a
classical result that 
\eqn\hard{Z_{{\rm hard\ dimers}\atop {\rm on\ [0,}n{\rm ]}}
= {1\over (1+y)^{n+1}}{1-y^{n+2}\over 1-y}}
with the parametrization
\eqn\Wy{W=-{1\over y+{1\over y}+2}\ .}
Since $W=-g_3^2\, R^3$, $y$ is precisely the solution of \ytriang.
\fig{Correspondence between configurations of non-intersecting paths
on the graph of Fig.~\triang\ where the last exit point is shifted by
one unit to the right, and hard dimers on a segment. We distinguish
between the cases where the uppermost point $(1,n+1)$ is attained
(bottom) or not (top).}{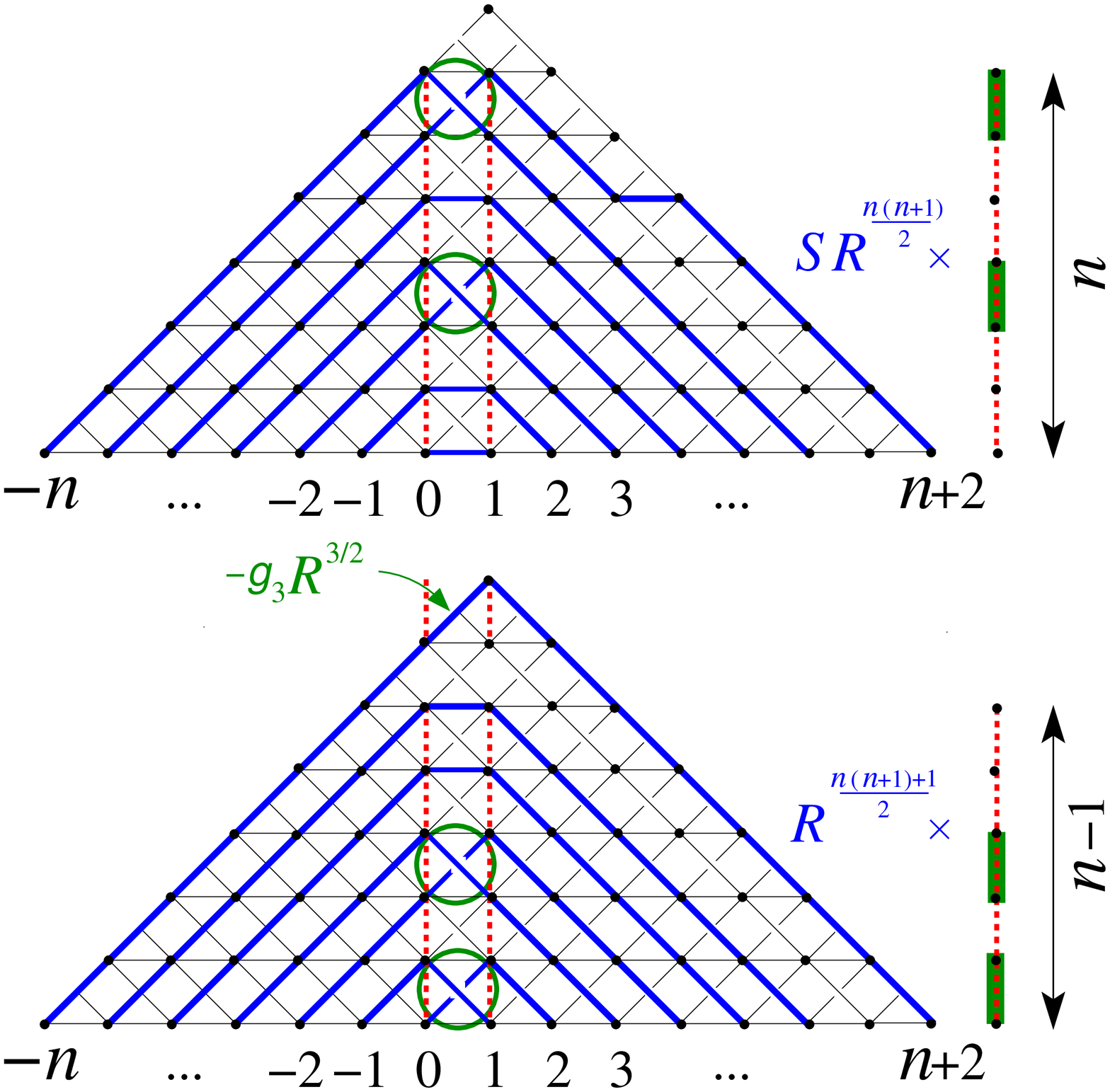}{10.cm} \figlabel\triangHOS
As for the Hankel minor $\tilde{H}_n$, it is easily identified in the 
hard dimer language as 
(see Fig.~\triangHOS\ for an illustration)
\eqn\Hntildetriang{\tilde{H}_n=(n+1)S\, R^{{n(n+1)\over 2}} 
Z_{{\rm hard\ dimers}\atop {\rm on\ [0,}n{\rm ]}}
-g_3 R^{3/2}\, R^{{n^2+n+1\over 2}}
Z_{{\rm hard\ dimers}\atop {\rm on\ [0,}n-1{\rm ]}}\ .}
Expressions \RnSntrian\ follow immediately via \Rndet\ and \Sndet.

\subsec{Quadrangulations}

\fig{Interpretation of formula (6.20) for $F_{2i+2j}$ in the case
  of quadrangulations. The weights per step are different (and allow
  for a level-step) in the central strip as
  shown.}{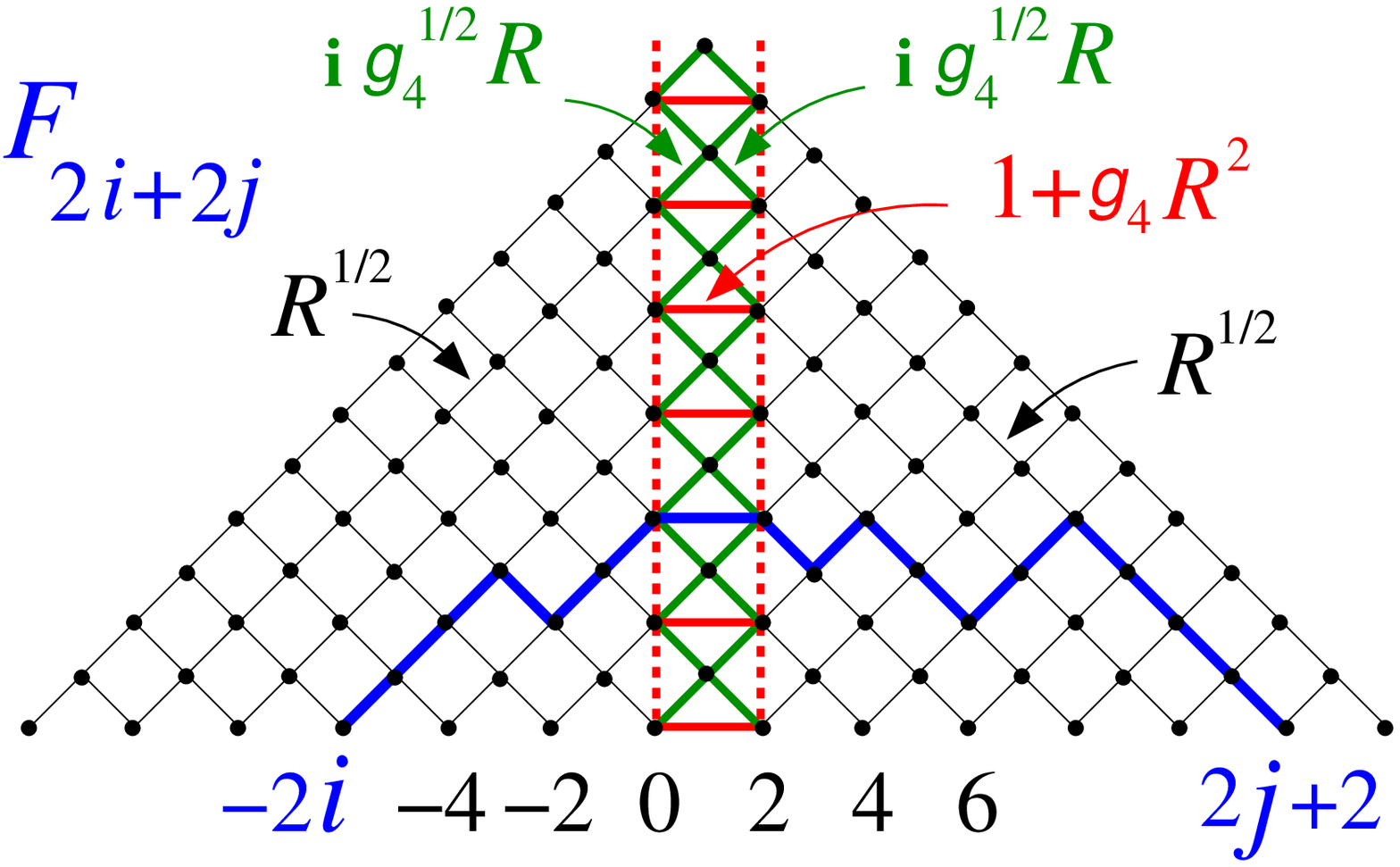}{8.cm} \figlabel\quadrang
\fig{Correspondence between configurations of non-intersecting paths
on the graph of Fig.~\quadrang\ and hard dimers on a segment of length
$2n+1$.}{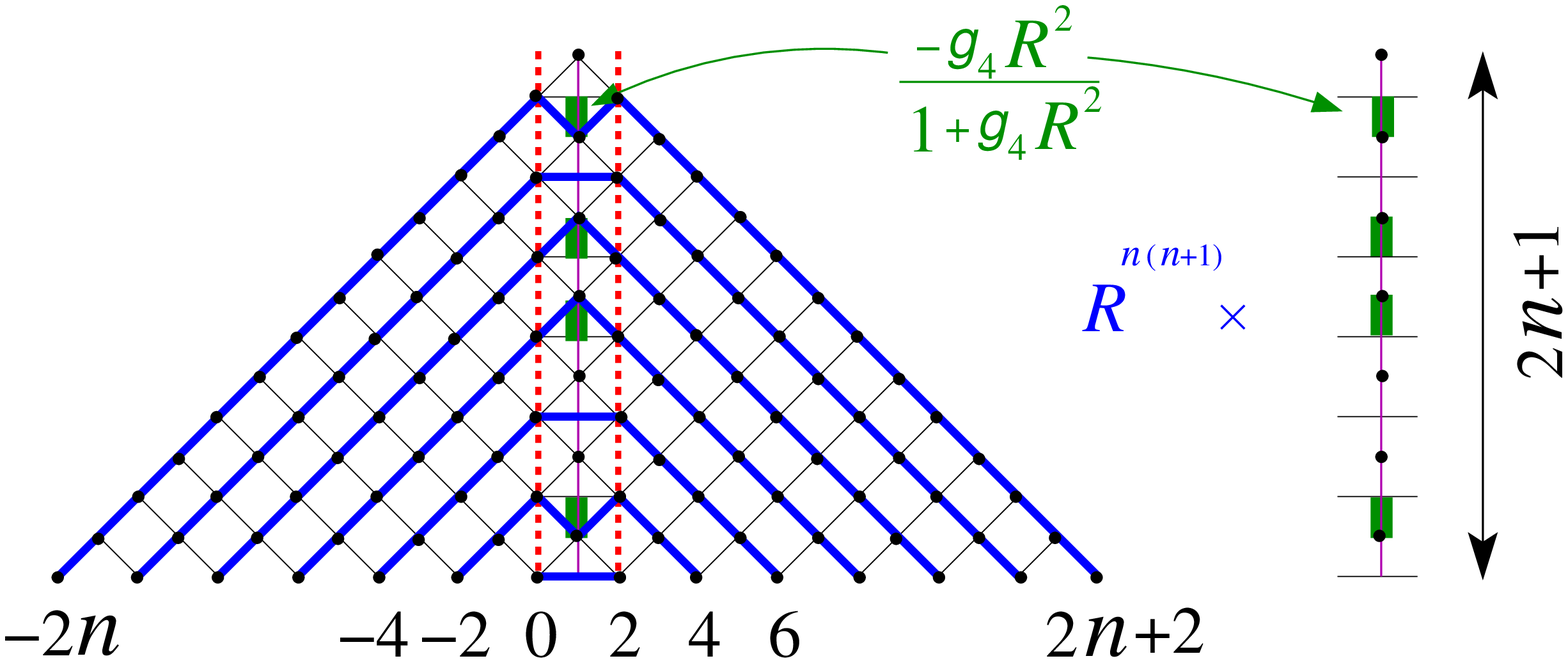}{12.cm} \figlabel\quadrangHO
Let us now come to the case of quadrangulations, characterized by 
$g_k= g_4 \delta_{k,4}$. We have $S=0$ in this case so that only
paths of even length contribute to $P$ or $P^+$. The only non-vanishing 
$A_q$ coefficients are $A_0$ and $A_2$, from \apaths, given by
\eqn\azerotwo{A_0=R-2 g_4\, R^2= 1+g_4\, R^2, \qquad A_2=-g_4\, R}
with $R$ now given, from \RSeqs, by 
\eqn\Rquad{R=1+3 g_4\, R^2\ .}
As quadrangulations belong to the class of maps with even face degrees, 
we use the results of Sect.~5. The non-vanishing ${\hat B}_i$'s are
\eqn\hatbquad{{\hat B}_0=1-g_4\, R^2,\quad {\hat B}_1={\hat B}_{-1}=-g_4\, R^2}
so that the characteristic equation reads
\eqn\charcheckquad{(1-g_4\, R^2)-g_4\, R^2\left(y+{1\over y}\right)=0\ .}
By \finalRnpair, the two-point function reads 
\eqn\Rnquadrang{R_n= R\, {(1-y^n)(1-y^{n+3})\over 
(1-y^{n+1})(1-y^{n+2})}\ ,}
which matches Eq.~(4.10) of \GEOD\ (up to the change of notation $y\to x$
and $n\to n+1$).

Again, we may give a combinatorial hard-dimer interpretation to
these results as follows.  
The generating function $F_n$ vanishes for $n$ odd and we have
\eqn\Ftwon{F_{2n}=A_0\, P^+(2n;R,0)+A_2\, P^+(2n+2;R,0)=A_0\, R^n {\rm cat_n}
+A_2\, R^{n+1} {\rm cat_{n+1}}}
where ${\rm cat}_n={2n \choose n}/(n+1)$ are the Catalan numbers.
In particular $F_{2i+2j}$ may be interpreted as enumerating paths from, 
say $(-2i,0)$ to $(2j+2,0)$ on the graph made of the restriction
in the upper half-plane made of a tilted square grid, completed
by horizontal segments in the central strip between abscissas $0$ and $2$
(see Fig.~\quadrang).
Again the graph is implicitly oriented from left to right so that paths
have increasing abscissas. Each up- or down-step of the path receives
a factor $\sqrt{R}$ except those in the central strip which receive
instead a weight $\sqrt{R}\sqrt{A_2}= {\rm i} \sqrt{g_4} R$, while the
horizontal paths receive a weight $A_0=1+g_4 R^2$. Using again the 
LGV lemma \LGV,
the quantity $h^{(0)}_n=\det_{0\leq i,j\leq n} F_{2i+2j}$ enumerates sets
of $n+1$ non-intersecting paths from points $(I_j)_{j=0,\cdots n}$ with 
coordinates $(-2j,0)$ to points $(O_j)_{j=0,\cdots n}$ with coordinates
$(2j+2,0)$. Again these paths have fixed ascending and descending parts
on both sides of the central strip, contributing an overall factor 
$R^{n(n+1)}$ and the only freedom comes from the central strip 
where the $j$-th path (numbered $0$ to $n$ from bottom to top) 
connects $(0,2j)$ to $(2,2j)$ either via a horizontal step with
weight $1+g_4 R^2$, or by a two-step sequence passing either by $(1,2j-1)$ (if 
$j\geq 1$) or by $(1,2j-1)$, with a total weight $({\rm i} \sqrt{g_4} R)^2=
-g_4 R^2$. Since the paths are non-intersecting, these (up or down) two-step
sequences cannot be adjacent in the strip and their vertical positions
define a hard dimer configuration in $[0,2n+1]$ (see Fig.~\quadrangHO). 
Extracting an overall factor $(1+g_4 R^2)^{n+1}$, each dimer receives a weight
\eqn\Wval{W=-{g_4 R^2\over 1+g_4 R^2}\ .}
so that we have 
\eqn\hnzero{h_n^{(0)}= R^{n(n+1)} (1+g_4 R^2)^{n+1} 
Z_{{\rm hard\ dimers}\atop {\rm on\ [0,}2n+1{\rm ]}}\ .}
%\eqn\valyquad{y+{1\over y}+1={1\over g_4 R^2}\ .}
\fig{Interpretation of formula (6.20) for $F_{2i+2j+2}$ in the
  case of quadrangulations. Note that the central strip is shifted by
  one unit with respect to Fig.~\quadrang}{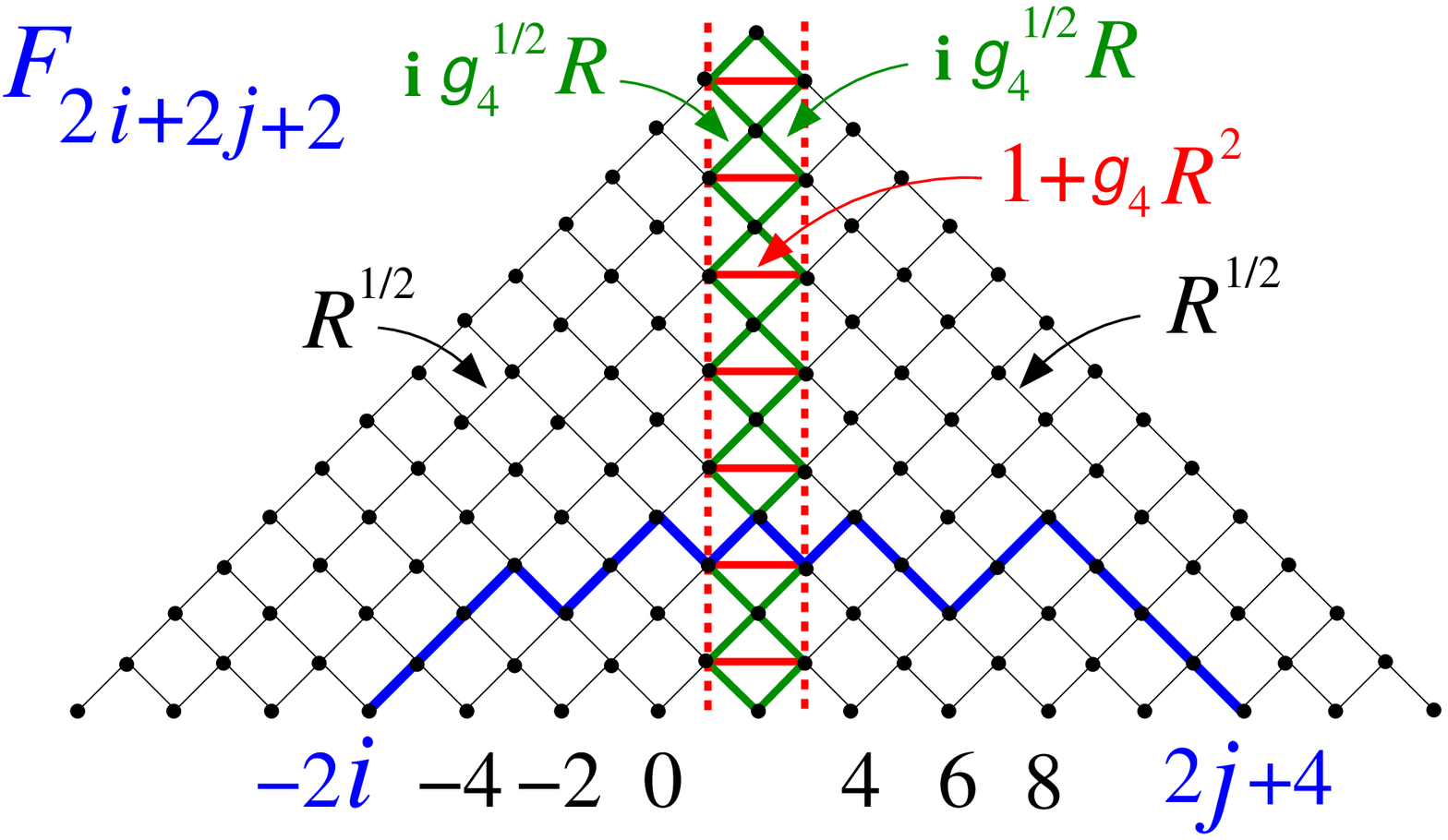}{8.cm}
  \figlabel\quadrangbis
\fig{Correspondence between configurations of non-intersecting paths
on the graph of Fig.~\quadrangbis\ and hard dimers on a segment of
length $2n+2$.}{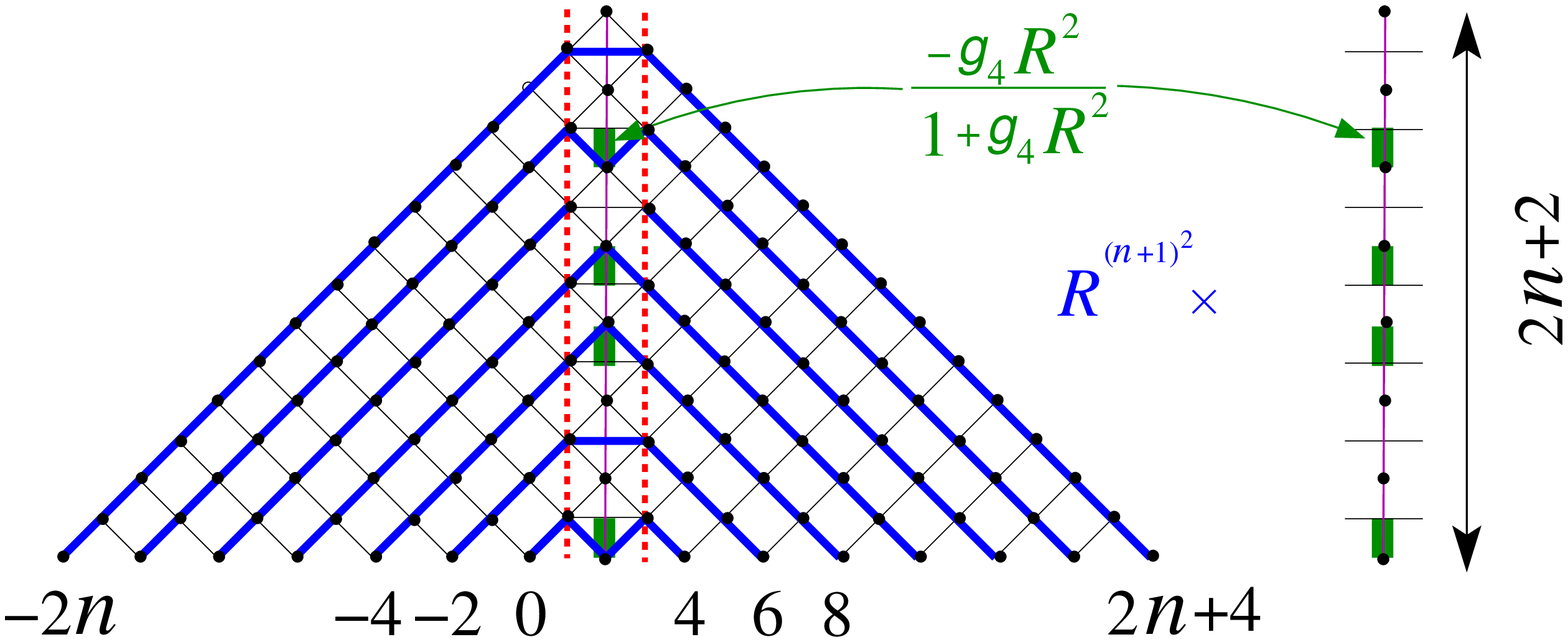}{12.cm} \figlabel\quadrangHObis
As for $h_n^{(1)}=\det_{0\leq i,j\leq n} F_{2i+2j+2}$, it enumerates sets
of $n+1$ non-intersecting paths from points $(I_j)_{j=0,\cdots n}$ with 
coordinates $(-2j,0)$ to points $(O_j)_{j=0,\cdots n}$ with coordinates
$(2j+4,0)$ with now a strip between abscissas $1$ and $3$ (see 
Fig.~\quadrangbis), easily transformed into a hard dimer configuration 
in $[0,2n+2]$ (see Fig.~\quadrangHObis). We deduce from this picture 
the formula
\eqn\hnone{h_n^{(1)}= R^{(n+1)^2} (1+g_4 R^2)^{n+1} 
Z_{{\rm hard\ dimers}\atop {\rm on\ [0,}2n+2{\rm ]}}\ .}
Using the explicit form \hard, this leads eventually to 
\eqn\Rquadrang{\eqalign{
R_{2n+1}&= {h_n^{(1)}\over h_{n-1}^{(1)}}/{h_n^{(0)}\over 
h_{n-1}^{(0)}}= R\, 
{(1-y^{2n+4})(1-y^{2n+1})\over 
(1-y^{2n+2})(1-y^{2n+3})} \cr
R_{2n}&= {h_n^{(0)}\over h_{n-1}^{(0)}}/{h_{n-1}^{(1)}\over 
h_{n-2}^{(1)}}= R\,  
{(1-y^{2n+3})(1-y^{2n})\over 
(1-y^{2n+1})(1-y^{2n+2})} \cr}
}
where $y$ is precisely the solution of \charcheckquad.
The above equations are summarized into \Rnquadrang.

\newsec{Conclusion and discussion}

To conclude, we discuss further connections of our results raising 
some open questions. Sect.~7.1
is devoted to the manifestations in our context of the intimate 
relation between continued fractions and orthogonal polynomials.
Sect.~7.2 makes the connection with the matrix integral approach.

\subsec{Convergents and orthogonal polynomials}

There is a deep connection between continued fractions and orthogonal
polynomials. Let us indeed consider the family of polynomials
$(q_n(z))_{n \geq 0}$ defined by
\eqn\expliorth{q_n(z) = {1 \over H_{n-1}} \det_{0 \leq i,j \leq n}
 \left( F_{i+j} \, \raise -2pt \vdots\,  z^j \right)}
where $(a_{i,j}\, \raise -2pt \vdots\, a_j)$ denotes the matrix with
elements $a_j$ in the last row and elements $a_{i,j}$ in all previous
rows.  It is a classical exercise to check that $q_n(z)$ is a monic
polynomial of degree $n$, that the family is orthogonal with
respect to the scalar product defined by $\langle z^n,z^m
\rangle=F_{n+m}$, with $\langle q_n(z),q_m(z)
\rangle= \delta_{n,m} H_n / H_{n-1}$, and that it satisfies
Favard's three-term recurrence
\eqn\favard{z\, q_n(z) = q_{n+1}(z) + S_n\, q_n(z) + R_n\, q_{n-1}(z)}
with the initial data $q_{0}(z)=1$, $q_{-1}(z)=0$. The connection with
continued fractions comes from the fact that the reciprocal of
$q_n(z)$ appears as the denominator of the $n$th convergent of the
J-fraction \continued, namely
\eqn\convergent{{\tilde{p}_n(z) \over \tilde{q}_n(z)} =
  {1 \over {\displaystyle 1 - S_0 z - {R_1 z^2
  \over {\displaystyle \ddots \over
   \displaystyle 1 - S_{n-1} z - {R_n z^2 \over 1 - S_n z} }}}}}
where $\tilde{q}_n(z)=z^n q_n(1/z)$ and  where $\tilde{p}_n(z)$ are the
so-called numerator polynomials. $\tilde{p}_n(z)$ and $\tilde{q}_n(z)$
also appear within an expression for the $n$th truncation \truncation,
see for instance [\xref\FLAJFRAC,\xref\FLAJFRACBIS].

In the context of maps, the $d$th convergent of \continued\ may be
interpreted as the generating functions for rooted maps where every
vertex incident to the root face is at a distance lesser than or equal
to $d$ from the origin of the root edge. As for the $d$th truncation
\truncation, its interpretation is given in Sect.~3.3 via \Fnd\ as
a generating function for pointed rooted maps with a control both
on the root degree and on the distance from the origin to the root face.
Truncations are known to admit a simple expression in the case of
quadrangulations \PSEUDOQUAD, namely
\eqn\truncquad{{1 \over {\displaystyle 1 - {R_{d+1} z^2
\over \displaystyle 1 - {R_{d+2} z^2 \over 1 - \cdots}}}}= 
{\tilde W} {1-({\tilde W}-1)y{1-y^{d+1}\over
1-y^{d+3}}\over 1-({\tilde W}-1)y{1-y^{d}\over1-y^{d+2}}}}
with $R$ and $y$ as in Sect.~6.2 and
${\tilde W}=(1-\sqrt{1-4Rz^2})/(2 R z^2)$.
We may wonder if this expression
generalizes to arbitrary $g_k$'s. We provide in Appendix C a general
expression for the orthogonal polynomials $q_n(z)$. Unfortunately we lack an
expression for the $\tilde{p}_n(z)$, as they do not admit a formula
similar to \expliorth.

\subsec{Connection with matrix integrals}

Another classical and fruitful approach to map enumeration problems is
via matrix integrals. Let us now comment informally on its connection
with our present results.

Matrix integrals give a simple expression for the all-genus generating function
\eqn\FnN{F_n[N]\equiv \sum_{h\geq 0} N^{-2h}\, F_n^{(h)}}
where $F_n^{(h)}$ denotes the generating function for rooted maps of genus $h$
with a root face of degree $n$, with face weights $(g_k)_{k\geq 1}$ as before. 
In particular, we have $F_n=F_n^{(0)}= \lim_{N\to \infty} F_n[N]$. We
have the matrix integral representation
\eqn\FnNmatrix{F_n[N]= 
{\int\, dM\, \Tr\, (M^n) \exp(-N \Tr\, V(M))\over
N\ \int\, dM\, \exp(-N \Tr\, V(M))}}
where $dM$ denotes the Lebesgue (translation-invariant) measure over
the space of $N\times N$ hermitian matrices and
\eqn\defV{V(x)\equiv {x^2\over 2}-\sum_{k\geq 1} g_k {x^k\over k}\ .}
The terms $g_k x^k$ act as a ``perturbation'' of the quadratic
potential $x^2/2$ corresponding to the well-known Gaussian Unitary
Ensemble. In all rigor, expression \FnNmatrix\ must be understood as a
power series in the $g_k$'s and $N$, whose precise definition is beyond
the scope of this section.

We now briefly discuss the usual approaches for studying \FnNmatrix,
beside the loop equations already mentioned in Sect.~3.2. The
original ``physical'' approach is the so-called saddle-point or
steepest descent method \BIPZ. It consists in remarking that the integrands
in \FnNmatrix\ only depend on the eigenvalues of $M$, and observing
that for large $N$ the dominant contribution comes the ``equilibrium''
continuous distribution of eigenvalues. It yields \eqn\Fneigendens{F_n
= \int d\lambda \, \lambda^n \rho(\lambda)} where $\rho(\lambda)$
denotes the density of eigenvalues. Therefore, in the saddle-point
picture our moments $F_n$ are precisely those associated with the
spectral measure.

This has to be contrasted with the usual method of orthogonal
polynomials in random matrix theory \DGZ. There, we consider the family of
polynomials $(q_i^{(N)}(\lambda))_{i\geq 0}$,
with $q_i^{(N)}(\lambda)$ monic of degree
$i$, orthogonal with respect to the scalar product defined by
$(\lambda^n,\lambda^m)=\int d\lambda\, \lambda^{n+m} \exp(-N\, V(\lambda))$.
In other words, these orthogonal polynomials are defined with respect to
the ``$N=1$'' eigenvalue density while those of \expliorth\ are defined with 
respect to the ``$N \to \infty$'' eigenvalue density.
Favard's theorem states that we still have a three-term recurrence
\eqn\relpol{\lambda\, q_i^{(N)}(\lambda)= q_{i+1}^{(N)}(\lambda)
+S_i^{(N)}\, q_i^{(N)}(\lambda)+R_i^{(N)}\, q_{i-1}(\lambda)}
where $(S_i^{(N)})_{i\geq 0}$ and $(R_i^{(N)})_{i\geq 1}$ (with
$R_0^{(N)}=0$) are associated with the scalar product at hand, for
instance they may be expressed via Hankel determinants of the moments
of the measure $d\lambda\, \exp(-N\, V(\lambda))$. Now, because of the 
specific form of the scalar product $(\cdot,\cdot)$, 
we also have the relations \DGZ
\eqn\risiN{\eqalign{R_i^{(N)}& ={i\over N} + \sum_{k\geq 2} g_k\,
Z_{i,i-1}^{(N)}(k-1) \cr S_i^{(N)}& =\sum_{k\geq 1} g_k\,
Z_{i,i}^{(N)}(k-1)\ .\cr}}
Here $Z_{i,j}^{(N)}(k)$ denotes the
generating function of three-step paths from $(0,i)$ to $(k,j)$ with a
weight $R_{m}^{(N)}$ per down-step $(t,m)\to (t+1,m-1)$ and a weight
$S_m^{(N)}$ per level-step $(t,m)\to (t+1,m)$.  Remarkably enough,
these equations look very similar to \rnsnrecur, and
allow to identify $R_i^{(N)}$ and $S_i^{(N)}$ with
generating functions for mobiles (respectively half-mobiles) as
defined in this paper with, however, a somewhat mysterious and
non-conventional weight $(m/N)$ for labeled vertices with label
$m$.
In terms of orthogonal polynomials, $F_n[N]$ is given by
\eqn\Fnortho{F_n[N] = {1\over N}\sum_{i=0}^{N-1} { \left( \lambda^n
  q_i^{(N)}(\lambda), q_i^{(N)}(\lambda) \right) \over \left(
  q_i^{(N)}(\lambda), q_i^{(N)}(\lambda) \right) } =
  {1\over N} \sum_{i=0}^{N-1} Z_{i,i}^{(N)}(n)\ . }

We may further connect the relations \risiN\ to known 
mobile generating functions by
considering the limit $N\to \infty$. In this limit, the quantity
$u=i/N$ is treated as a continuous variable, and $R_i^{(N)}$ and
$S_n^{(N)}$ become at leading order in $1/N$ smooth functions $R(u)$
and $S(u)$. Then, in this limit \risiN\ yields
\eqn\impliru{\eqalign{R(u)&= u +\sum_{k\geq 2} g_k \sqrt{R(u)} P_{-1}(k-1;
R(u),S(u)) \cr
S(u)&= \sum_{k\geq 1} g_k P(k-1; R(u),S(u)) \cr}}
involving paths which receive homogeneous weights $R(u)$ and $S(u)$,
irrespectively of the ordinates. Setting $u=1$, we recover precisely
\rsrecurbis. For general $u$, $R(u)$ and $S(u)$ are generating functions for
 mobiles and half-mobiles with an extra weight $u$ per labeled vertex,
as encountered in Sect.~3.1 and Appendix B.
In the limit $N\to \infty$, the sum in \Fnortho\ becomes an integral, 
leading to the expression
\eqn\Fnortholimit{F_n= \int_0^1 du\, P(n;R(u),S(u))} 
which is consistent with the discussion of Sect.~3.1.

For finite $N$, a natural question is whether the unconventional weight 
of labeled vertices, as well as the form \Fnortho\ of the all-genus 
generating function may be given a direct combinatorial interpretation.
This could open the way to get explicit formulas for discrete 
distance-dependent two-point functions in maps of higher genus.

\bigskip
\noindent{\bf Acknowledgments:}
We thank C. Krattenthaler for pointing out the connection with Schur
functions.  Part of this work was completed at the Centre \'Emile Borel of
Institut Henri Poincar\'e.

\appendix{A}{Decomposition of maps into slices}

\fig{The decomposition (a) of a pointed rooted map into slices by
cutting it along all leftmost geodesic paths emerging from vertices
incident to the root face. Slices are of two types: those (b)
with a basal edge of type $(m,m)$, counted by $S_m$, and those (c)
with a basal edge of type $(m,m-1)$, counted by $R_m$. Note that
edges of type $(m-1,m)$ delimit empty slices. The actual depth of a
slice may be less than $m$ since its two boundaries may merge before
the origin ($s_i$ and $r_i$ counting the slices with depth
$i$).}{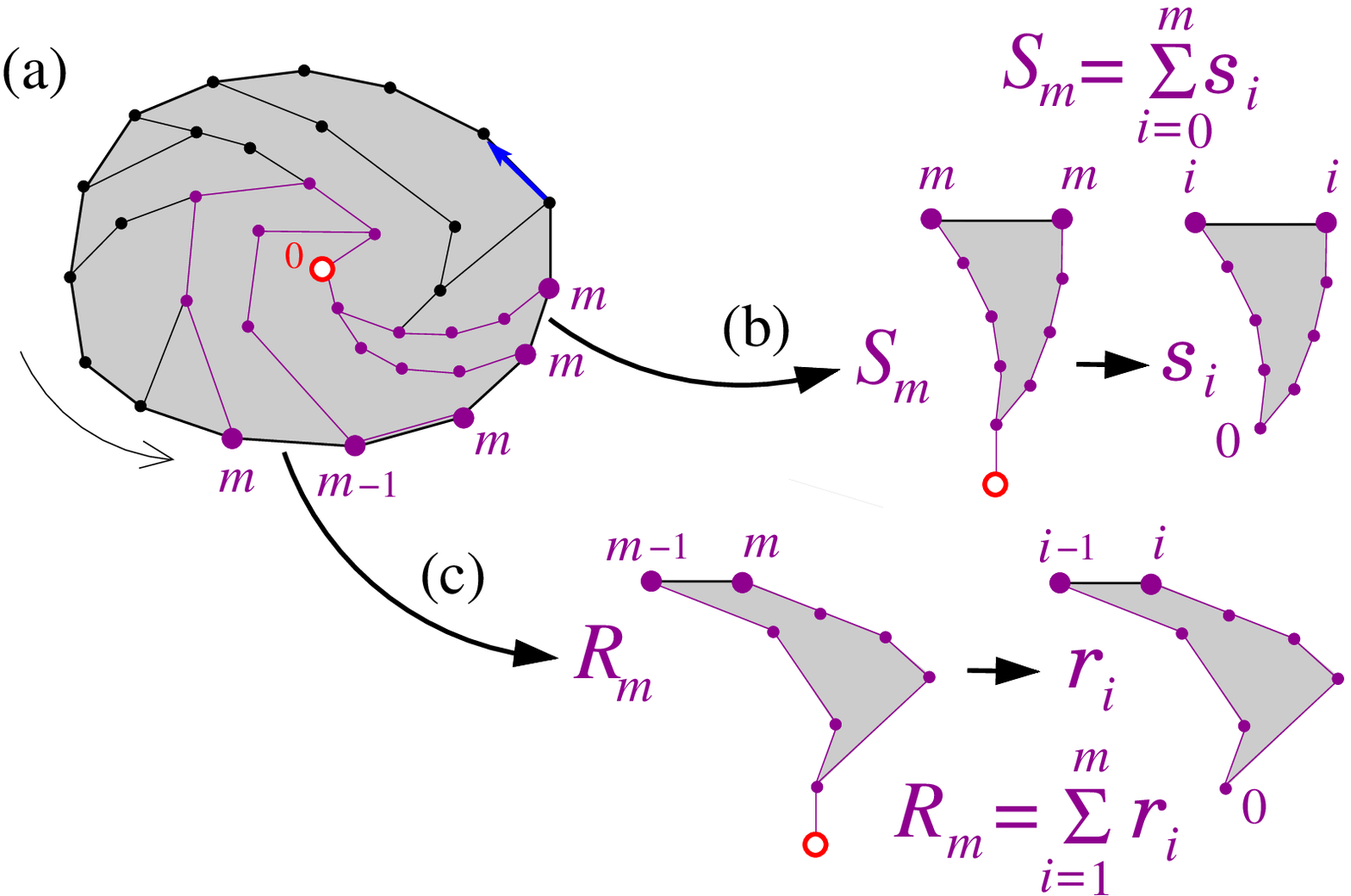}{12.cm} \figlabel\slicester

So far we used mobiles as the natural framework to introduce continued
fractions in map enumeration problems. Still, the explicit recourse to
mobiles is not {\it stricto sensu} necessary to prove the basic
continued fraction expansion \continued\ and more generally the
interpretation \Fnd\ of truncations \truncation, when the weights
$R_m$ and $S_m$ have their original interpretation of Sect.~1 as
map generating functions. This relation can be understood
alternatively as the result of some appropriate decomposition of maps
which can be defined without reference to the underlying mobile
structure even though, in practice, the two are intimately related.

We start again with pointed rooted maps enumerated by $f_{n;d}$, i.e. maps
with a root face of length $n$, with the origin of the root edge at distance 
$d$ from the origin of the map and with all the other vertices incident to the 
root face at a distance larger than or equal to $d$. As before, we choose the 
root face for external face and the clockwise sequence of distances of
its incident vertices defines a three-step path of length $n$ with endpoints
at ordinate $d$, and which stays above $d$. For each vertex incident to the 
root face, we draw the {\it leftmost geodesic} path from this vertex to the 
origin of the map. Cutting along these paths decomposes the map into a 
number of connected domains that we call {\it slices} 
(see Fig.~\slicester). These slices are maps whose boundary is made of one
{\it basal edge} originally incident to the root face and of two leftmost 
geodesic 
paths merging at some apex (possibly different from the former origin of the 
map). Note that the leftmost geodesics may themselves contain edges 
originally incident to the root face. The basal edge of 
a non-empty slice is an edge originally of type $(m,m)$ or of type $(m,m-1)$ 
(counterclockwise). Indeed, the edges of type $(m-1,m)$ give rive to empty 
slices as they lie on the leftmost geodesic path emerging from their endpoint 
at distance $m$. Non-empty slices are therefore associated only with the down- and 
level-steps of the three-step path above. Since two consecutive leftmost 
geodesics may merge before reaching the origin of the map, the actual 
length of the leftmost geodesics is reduced in the slice by some value $m-i$ 
corresponding to the length of their common part. This leads us to associate 
a weight $S_m=\sum_{i=0}^m s_i$ (respectively $R_m=\sum_{i=1}^m r_i$) 
to each level- (respectively down-) step starting at ordinate $m$ in 
the three-step path, where $s_i$ (respectively $r_i$) are the generating 
functions for slices of depth $i$, i.e. with a boundary made of a basal edge 
and two leftmost 
geodesics of the same length $i$ (respectively of length $i-1$ and $i$)
(see Fig.~\slicester). Here again, we have to add to $r_1$, hence to
$R_m$, a conventional factor $1$ which accounts for the case where an edge 
of type $(m,m-1)$ incident to the root face would be the boundary of an
empty slice, which happens when its endpoint at distance $m-1$ is 
the only vertex 
at distance $m-1$ from the origin adjacent to its endpoint at distance $m$.
On the contrary, each edge of type $(m,m)$ incident to the root 
face is the basal edge of a non-empty slice. We therefore recover precisely the 
generating function $Z^+_{d,d}(n)$ defined in Sect.~2.1, which, from the above 
decomposition procedure, enumerates all maps in $f_{n;d}$, as well as all 
maps in $f_{n;i}$ for $i<d$, since the concatenation of slices produces maps 
whose origin may remain at a distance less than $d$. We end up with 
the desired relation $F_{n;d}= Z^+_{d,d}(n)$, but now with a different 
interpretation for the weights $S_m$ and $R_m$ as slice generating functions. 
It remains to show that these new definitions match their former definition
of Sect.~1.2 in terms of map generating functions.

\fig{Gluing two slices in $s_i$ and $s_j$ and identifying pairwise
boundary edges creates a pointed rooted map of type $(m \to m)$ with
$m = \max(i,j)$, implying $S_n^2=T_n$ with $T_n$ defined as in Sect.~1.
Identifying pairwise boundary edges of a slice in $r_i$ creates a
pointed rooted map of type $(i \to i-1)$.}{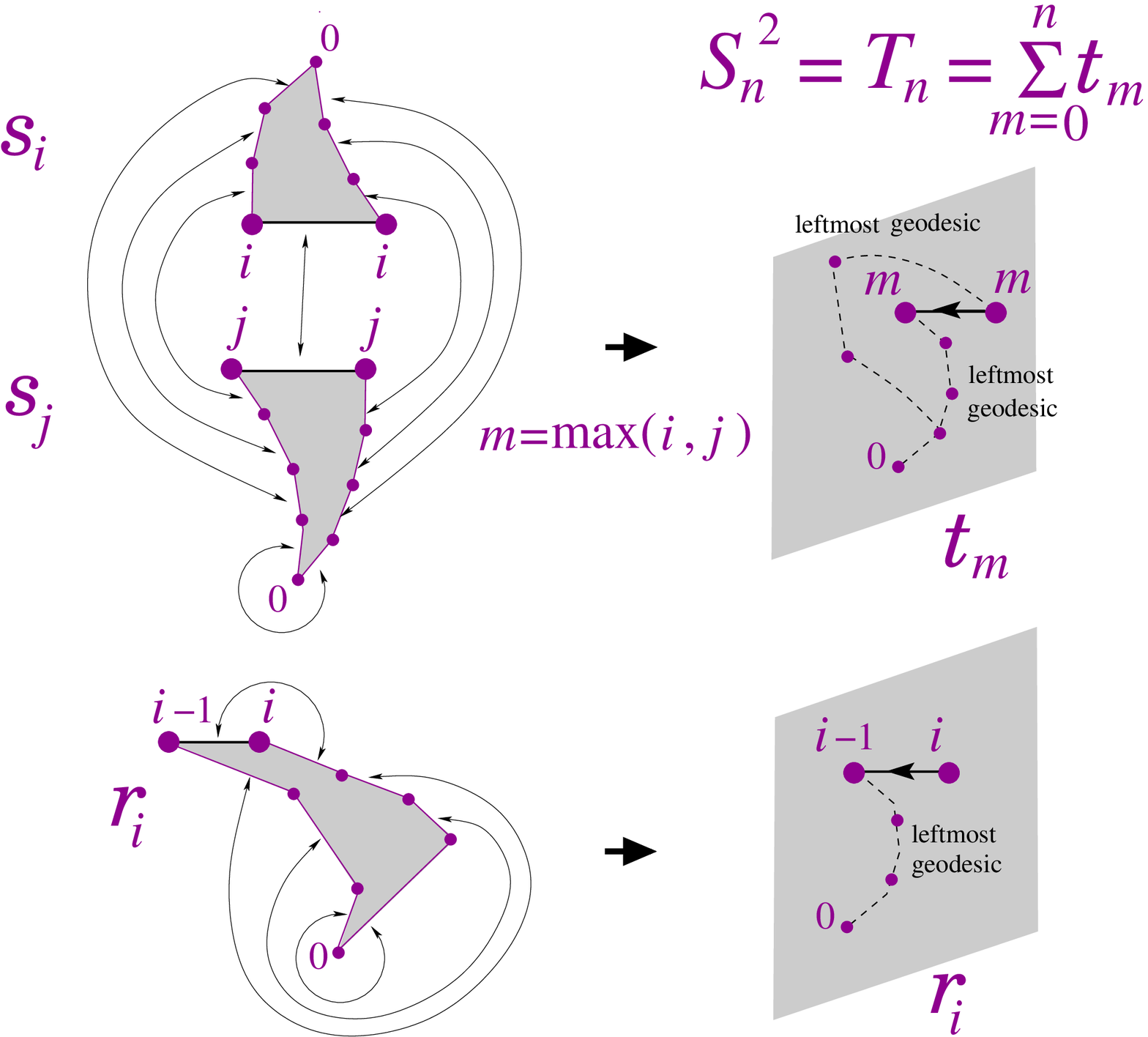}{12.cm}
\figlabel\slicestwoter
It is easily seen (Fig.~\slicestwoter) that, starting from a slice counted  
by $r_i$, completing the leftmost geodesic of length $i-1$ by the 
basal edge (of type $(i-1,i)$) and gluing it with the other leftmost 
geodesic of length $i$ produces bijectively a pointed map with a marked 
edge of type $(i,i-1)$ with respect to the origin of the map. 
This marked edge is promoted as root edge.
Note that, after gluing, we are apparently left with a marked geodesic
but it is by construction the leftmost geodesic starting with the
root edge (oriented from $i$ to $i-1$) and ending
at the origin. It is therefore uniquely determined
and may be erased without loss of information. This bijection guarantees
that our new definition of $r_i$ (and $R_m$) matches the definition 
of Sect.~1. 
As for $S_n$, we see (Fig.~\slicestwoter) that, upon gluing two 
slices in $S_n$ (respectively in $s_i$ and $s_j$ for some 
$i\leq n$ and $j\leq n$), we get a pointed map with a marked
edge of type $(m,m)$ for some $m\leq n$ (with $m=\max(i,j)$).
This marked edge is oriented with, say, the first slice on its right, and
promoted as root edge. 
Again, we have after gluing two marked leftmost 
geodesics (starting respectively with the root edge and the reversed
root edge) which may be erased without loss of information. 
Hence $S_n^2=T_n$ with $T_n$ defined as in Sect.~1.

\appendix{B}{Proof of \Fnu}

Let us start by showing the following identity, valid for any $n\geq 0$ and
$k\geq 2$
\eqn\identun{\eqalign{{\partial\ \over \partial R}\left(R \sum_{q=0}^{k-2}
P^+(n+q;R,S)P(k-2-q;R,S)\right) & = P(n;R,S) {\partial\ \over \partial R}\left(
\sqrt{R}\, P_{-1}(k-1;R,S)\right)\cr &\ \ \ + \sqrt{R}\, P_{-1}(n;R,S)
{\partial\ \over \partial R} P(k-1;R,S)\ .\cr}}
At this stage $R$ and $S$ are arbitrary coefficients. Attaching a weight
$x^n y^{k-2}$ and summing over $n\geq 0$ and $k\geq 2$, the 
r.h.s produces the combination
\eqn\rhs{\pi(x)\ {\partial\ \over \partial R} \left({\pi_{-1}(y)\over y}\right)+
\pi_{-1}(x)\ {\partial\ \over \partial R}\left({\pi(y)-1\over y}\right)\ ,}
where we introduced the generating functions
\eqn\pione{\eqalign{\pi(z)& \equiv 
\sum_{n\geq 0}P(n;R,S)z^n= {1\over \sqrt{\kappa(z)}}
\cr \pi_{-1}(z)& \equiv \sum_{n\geq 0}\sqrt{R}\, P_{-1}(n;R,S)z^n
= z {1-Sz -\sqrt{\kappa(z)} \over 2 Rz^2\sqrt{\kappa(z)}} \cr}}
with $\kappa(z)$ as in \kappasr. The last formula comes from the identification
$\pi_{-1}(z)=z \pi^+(z)\pi(z)$, where
\eqn\piplus{\pi^+(z)\equiv \sum_{n\geq 0} P^+(n;R,S) z^n =
{1-Sz-\sqrt{\kappa(z)} \over 2 Rz^2}\ .}
As for the l.h.s in \identun, setting $s=n+q$ and $t=k-2-q$,
we have to compute
\eqn\changsum{\eqalign{
\sum_{n\geq 0}\sum_{k\geq 2}\sum_{q=0}^{k-2} P^+(n+q;R,S) & P(k-2-q;R,S) 
x^n y^{k-2} \cr 
& = \sum_{s\geq 0} \sum_{t\geq 0}\sum_{q=0}^s x^s y^t 
\left({y\over x}\right)^q
P^+(s;R,S) P(t;R,S) \cr  
& = \sum_{s\geq 0}\sum_{t\geq 0} {x^{s+1}-y^{s+1}\over x-y}\, y^t
P^+(s;R,S) P(t;R,S) \cr
& = {x \pi^+(x)- y \pi^+(y) \over x-y}\, \pi(y)\ .\cr}}
Proving \identun\ therefore reduces to checking the relation
\eqn\relun{{\partial\ \over \partial R}\left(R\, {x \pi^+(x)- y \pi^+(y) \over x-y}\, \pi(y)
\right) = 
\pi(x)\ {\partial\ \over \partial R} \left({\pi_{-1}(y)\over y}\right)+
\pi_{-1}(x)\ {\partial\ \over \partial R}\left({\pi(y)-1\over y}\right)\ ,}
which, knowing the explicit forms of all the involved generating functions, 
is a straightforward task.  

Multiplying \identun\ by $\delta_{k,2} - g_k$, summing over $k\geq 2$ and
exchanging the sums over $k$ and $q$ leads to 
\eqn\dfnR{{\partial\ \over \partial R} \sum_{q=0}^{\infty} A_q P^+(n+q;R,S)
   = P(n;R,S) {\partial\ \over \partial R} u + \sqrt{R}\, 
  P_{-1}(n;R,S) {\partial\ \over \partial R}v}
where we introduced
\eqn\defuv{\eqalign{ 
u & \equiv R- \sum_{k\geq 2} g_k\, \sqrt{R}\, P_{-1}(k-1;R,S) \cr
v & \equiv S- \sum_{k\geq 1} g_k\, P(k-1;R,S)
\cr }}
and where $A_q$ is defined as in \apaths. Note that we added 
for convenience a trivial constant term $g_1$ in the definition of $v$,
which disappears after derivation with respect to $R$. 

Similarly, one can easily prove 
\eqn\identdeux{\eqalign{{\partial\ \over \partial S}\left(R \sum_{q=0}^{k-2}
P^+(n+q;R,S)P(k-2-q;R,S)\right) & = P(n;R,S) {\partial\ \over \partial S}\left(
\sqrt{R}\, P_{-1}(k-1;R,S)\right)\cr &\ \ \ + \sqrt{R}\, P_{-1}(n;R,S)
{\partial\ \over \partial S} P(k-1;R,S)\cr}}
by checking the identity 
\eqn\reldeux{{\partial\ \over \partial S}\left(R\, {x \pi^+(x)- y \pi^+(y) \over x-y}\, \pi(y)
\right) = 
\pi(x)\ {\partial\ \over \partial S} \left({\pi_{-1}(y)\over y}\right)+
\pi_{-1}(x)\ {\partial\ \over \partial S}\left({\pi(y)-1\over y}\right)\ .}
This implies the relation
\eqn\dfnS{{\partial\ \over \partial S} \sum_{q=0}^{\infty} A_q P^+(n+q;R,S)
  = P(n;R,S) {\partial\ \over \partial S} u + \sqrt{R}\, 
  P_{-1}(n;R,S) {\partial\ \over \partial S}v}
with $u$ and $v$ as above. 
Equations \dfnR\ and \dfnS\ may be rewritten in terms of differential forms
\eqn\gradient{
d \left(\sum_{q=0}^{\infty} A_q P^+(n+q;R,S)\right)
= P(n;R,S)\, du +\sqrt{R} P_{-1}(n;R,S)\, dv\ .}
When $v=0$, Eqs.\defuv\ match precisely
the equations satisfied by the mobile generating functions $R(u)$ and
$S(u)$ with a weight $u$ per labeled vertex. 
Therefore, by restricting \gradient\ to the line $v=0$, we deduce
\eqn\dfnu{{d\ \over du} \sum_{q=0}^{\infty} A_q(u) P^+(n+q;R(u),S(u))
  = P(n;R(u),S(u))}
with $A_q(u)$ defined as in \apaths\ with $R,S$ replaced by $R(u),S(u)$.
This establishes Eq.~\Fnu\ upon integrating with respect to $u$ 
and noting that both sides of \Fnu\ vanish at $u=0$ ($A_q(0)=0$ 
since $R(0)=0$ and $F_n(0)=0$ from its map definition).

\appendix{C}{An expression for the orthogonal polynomials}

Using the notations of Sect.~4, we have $q_n(z) = {\displaystyle {1
\over H_{n-1}}} \det \pmatrix{ {\bf H}_n^{(n;)} \cr {\bf z}_n}$ where
${\bf z}_n$ is the row vector $\left( z^j \right)_{0 \leq j \leq n}$.
Inspired by \matrixid, we may write 
\eqn\orthmat{ 
  \pmatrix{ {\bf H}_n^{(n;)} \cr {\bf z}_n } =
  \pmatrix{ {\bf T}_{n-1}^{\rm T} & 0 \cr 0 & 1} \cdot
  \pmatrix{ {\bf B}_n^{(n;)} \cr {\bf z}_n \cdot {\bf T}_n^{-1} } \cdot
  {\bf T}_n.}
By Proposition 11 of \FLAJFRAC, the inverse matrix of ${\bf T}
=(P_i^+(j;R,S))_{i,j \geq 0}$ is given by the matrix of coefficients
of orthogonal polynomials associated with the continued fraction
\motzfrac. The $j$th such polynomial is $U_j((z-S)R^{-1/2})$, where
$U_j(z)$ denotes the $j$th Chebyshev polynomial of the second kind
(defined for instance via $U_j(2 \cos \theta) \sin \theta = \sin (j+1)
\theta$).  Hence we have ${\bf z}_n \cdot {\bf T}_n^{-1} = \left(
U_j((z-S)R^{-1/2}) \right)_{0 \leq j \leq n}$. Passing to determinants
in \orthmat, expanding the middle r.h.s determinant over the last row
and using the known expressions for $\det {\bf T}_n$ and $H_n$, we
obtain
\eqn\qnB{\eqalign{q_n(z) &= R^{n/2}
  \sum_{m=0}^n (-1)^{m+n} U_m\left({z-S \over \sqrt{R}}\right)
  { \det {\bf B}_n^{(n;m)} \over \det {\bf B}_{n-1}} \cr
  &= R^{n/2} \sum_{m=0}^n U_m\left({z-S \over \sqrt{R}}\right)
    {{\rm sp}_{2p}(\lambda_{p,n}^{(m)},{\bf x}) \over
     {\rm sp}_{2p}(\lambda_{p,n},{\bf x}) }}}
where $\lambda_{p,n}^{(m)}=n^{p-1}m$ denotes the partition with $p-1$ parts of
size $n$ and one part of size $m$. The $h$- or Weyl formulas allow to
further reexpress $q_n(z)$ in terms of $p \times p$ determinants,
where only one row depends on $z$.

\listrefs

\end